\documentclass[preprint,12pt,3p]{elsarticle}
\usepackage{setspace}
\usepackage{amsmath}
\usepackage{amssymb}
\usepackage{lineno}
\usepackage[dvipsnames]{xcolor}
\usepackage[utf8]{inputenc}
\usepackage{verbatim}
\usepackage[mathcal]{euscript}
\usepackage{graphicx}
\usepackage[english]{babel}
\usepackage{multicol}
\usepackage[utf8]{inputenc}  
\usepackage[T1]{fontenc}  
\usepackage{caption}
\usepackage[colorlinks=true ]{hyperref}
\usepackage[final]{pdfpages}
\usepackage{multicol}
\usepackage{float}
\usepackage{tikz}
\usepackage{tikz-3dplot} 
\usepackage{pdfpages}       
\usepackage{varioref} 
\usepackage{float}
 \usepackage{multicol}
\usepackage{multirow}
\usetikzlibrary{quotes,arrows.meta}
\usetikzlibrary{shapes,arrows}
\numberwithin{equation}{section}  
\usepackage{float}
\definecolor{camel}{rgb}{0.76, 0.6, 0.42}

\newtheorem{definition}{Definition}[section]

\newtheorem{Claim}[definition]{Claim}
\newtheorem{Remark}[definition]{Remark}
\newtheorem{example}[definition]{Example}

\newtheorem{proof}[definition]{Proof}
\newtheorem{Corollary}[definition]{Corollary}

\newcommand \bei {\begin{itemize}}
\newcommand \eei {\end{itemize}}
\newcommand \ubar u

\newcommand \del \partial

\newcommand \la \langle
\newcommand \ra \rangle 

\newcommand \auth    \textsc
\newcommand \be {\begin{equation}}
\newcommand \ee {\end{equation}}

\newcommand \bcor {\begin{Corollary}}
\newcommand \ecor {\end{Corollary}}

\newcommand \bpro {\begin{proof}}
\newcommand \epro {\end{proof}}
\newcommand \bdf {\begin{definition}}
\newcommand \edf {\end{definition}}
\newcommand \bex {\begin{example}}
\newcommand \eex {\end{example}}
\newcommand \bcl {\begin{Claim}}
\newcommand \ecl {\end{Claim}}
\newcommand \brm {\begin{Remark}}
\newcommand \erm {\end{Remark}}

\let\oldmarginpar\marginpar
\renewcommand\marginpar[1]{\-\oldmarginpar[\raggedleft\footnotesize #1]%
{\raggedright\footnotesize #1}}
\begin{document}
 \begin{frontmatter}



\title{Localized RBF methods for modeling infiltration using the Kirchhoff-transformed Richards equation}

\author[mpu]{Mohamed Boujoudar}
\author[mpu,uottawa]{Abdelaziz Beljadid\corref{mycorrespondingauthor}}
\cortext[mycorrespondingauthor]{Corresponding author}
\ead{abdelaziz.beljadid@um6p.ma}
\author[LMA]{Ahmed Taik}
\address[mpu]{Mohammed~VI Polytechnic University, Green City, Morocco}
\address[uottawa]{ University of Ottawa, Canada}
\address[LMA]{ University  Hassan~II, Casablanca, Morocco}

\begin{abstract} 
We develop a new approach to solve the nonlinear Richards equation based on the Kirchhoff transformation and localized radial basis function (LRBF) techniques. Our aim is to reduce the nonlinearity of the governing equation and apply LRBF methods for modeling unsaturated flow through heterogeneous soils. In our methodology, we propose special techniques which deal with the heterogeneity of the medium in order to apply the Kirchhoff transformation where we used the Brooks and Corey model for the capillary pressure function and a power-law relation in saturation for the relative permeability function. The new approach allows us to avoid the technical issues encountered in the Kirchhoff transformation due to soil heterogeneity in order to reduce the nonlinearity of the model equation. The resulting Kirchhoff-transformed Richards equation is solved using LRBF methods which have advantages in terms of computational cost since they don't require mesh generation. Furthermore, these LRBF techniques lead to a system with a sparse matrix which allows us to avoid ill-conditioned issues. To validate the developed approach for predicting  the dynamics of unsaturated flow in porous media, numerical experiments are performed in one, two, and three-dimensional soils. The numerical results demonstrate the efficiency and accuracy of the proposed techniques for modeling infiltration through heterogeneous soils. 
\end{abstract} 
\end{frontmatter}

\section{Introduction} 
Understanding infiltration through soils is of great importance in the fields of agriculture, hydrology, and water resources and environmental management. The modeling of infiltration processes is time-consuming and there is a need in the development of efficient techniques for these processes in the case of heterogeneous soils. The Richards equation~\cite{richards1931capillary} describes the dynamic of unsaturated flow through porous medium which is due to the actions of gravity and capillarity. Richards’ equation is highly nonlinear because of the largely nonlinear dependencies of both unsaturated hydraulic conductivity and capillary pressure on saturation \cite{gardner1958some, brooks1964hydraulic, van1980closed}. The van Genuchten \cite{van1980closed} and Brooks-Corey models \cite{brooks1964hydraulic} are often used for the capillary pressure function of unsaturated soils. In terms of numerical analysis, the Gardner model \cite{gardner1958some} is important since large class of analytical solutions are available \cite{SrivastavaYeh, tracy19951, huang2012analytical, hayek2016exact}. However, this model has some limitations in practical applications for describing unsaturated flow in soils \cite{rucker2005parameter}.

The design of efficient numerical techniques for solving the Richards equation is very challenging due to the highly non-linearity of the equation and the technical issues encountered in the numerical treatment of soil heterogeneity. Various classes of approaches have been developed for modeling infiltration in soils such as finite difference methods \cite{haverkamp1977comparison, celia1990general, clement1994physically}, finite element methods \cite{huyakorn1984techniques, radu2004order, bause2004computation, solin2011solving} and finite volume methods \cite{eymard1999finite, manzini2004mass,lai2015mass,svyatskiy2017second,ngo2020control}. While many methods have been developed for modeling unsaturated flow in soils, there is still a need for more efficient techniques to deal with the soil heterogeneity and the nonlinearity of the medium hydraulic properties as functions of saturation \cite{brooks1964hydraulic, van1980closed}. Most available approaches used iterative methods, such as Newton and Picard schemes \cite{celia1990general, kirkland1992algorithms, huang1996new, lehmann1998comparison,  an2012comparison, zha2017modified}, to linearize the system to be solved. While these iterative algorithms produce accurate results, they are expensive in terms of computational cost and they may have convergence issues for some flow conditions because of the highly non-linearity of the Richards equation and soil heterogeneity \cite{lehmann1998comparison,ji2008generalized,list2016study,zha2017modified}. 

Among the numerical techniques proposed to solve the Richards equation is the Kirchhoff transformation approach \cite{haverkamp1977comparison, pop2002error, radu2004order, ji2008generalized, berninger2011fast, suk2019numerical}. Several studies have shown the efficiency of this approach because it reduces the non-linearity of the Richards equation \cite{ross1990efficient, ji2008generalized, stevens2011scalable, boujoudar2022modelling}. However, this approach is developed for homogeneous soils or particular soil heterogeneity and capillary pressure functions \cite{haverkamp1977comparison, protopapas1991analytical, ji2008generalized, berninger2011fast, suk2019numerical}. Most of available techniques using the Kirchhoff transformation are limited to  Gardner model for the capillary pressure function \cite{yeh1989one,merrill1978laterally, tartakovsky2003unsaturated, bakker2004two, ji2008generalized, suk2019numerical, zhang2021finite}. For instance, Suk and Park \cite{suk2019numerical} recently developed a new numerical method based on the Kirchhoff transformation and the Gardner model to solve the Richards equation for layered soils. In their approach, the authors used a truncated Taylor series expansion to the Kirchhoff head at the material interface.

We note that the Gardner model is limited for practical applications compared to van Genuchten and Brooks-Corey models which are suitable for the entire range of pressure head \cite{suk2019numerical}. Incorporating both high non-linear models for the capillary pressure  \cite{brooks1964hydraulic,van1980closed} and soil heterogeneity is still a challenge in solving the Richards equation using the Kirchhoff transformation \cite{suk2019numerical}.

 This study is a follow-up of the paper \cite{boujoudar2021localized} in which we developed LRBF techniques for solving Richards equation in homogeneous medium where the Gardner model is used for capillary pressure. The techniques used in \cite{boujoudar2021localized} which allowed us to linearize the system can not be applied in the case of heterogeneous mediums with other models of capillary pressure such as the Brooks-Corey \cite{brooks1964hydraulic} model used in this study.
Here, we develop special techniques which deal with soil heterogeneity in order to apply the Kirchhoff transformation for solving the Richards equation based on the Brooks and Corey model \cite{brooks1964hydraulic} for capillary pressure. In our approach, a power-law relation in saturation for the relative permeability function is used to avoid the technical issues encountered in the Kirchhoff transformation due to the non-linearity of capillary pressure function and the heterogeneity of soils \cite{suk2019numerical}.

In our approach, we based on RBF meshless techniques \cite{kansa1990multiquadrics,lee2003local,li2013localized} to solve the resulting system. These techniques don't require mesh generation and are based only on a set of independent points, which makes them advantageous in terms of computational cost. Due to their simplicity to implement, they represent an attractive alternative to the classical methods as a solution method for partial differential equations. Note that there are two versions of RBF methods. The global method \cite{kansa1990multiquadrics1} and the local one \cite{lee2003local}. Several studies have demonstrated the efficiency of the global method \cite{kansa1990multiquadrics1,kansa1990multiquadrics1,mirinejad2017rbf} however it suffers from two major drawbacks: the ill-conditioned matrix obtained after the discretization process and the problem of choosing the adequate shape parameter for some RBFs \cite{vsarler2007global,li2013localized,hamaidi2016space}. To overcome these issues, the local methods were suggested \cite{lee2003local,li2013localized}. The LRBF methods have advantages in terms of operational memory and calculation procedures where only inversion of sparse matrix are required. On the other hand, the LRBF methods are efficient in solving high-dimensional problems with complex boundaries \cite{vsarler2007global,li2013localized,hamaidi2016space,boujoudar2021localized} and are less sensitive to the choice of the shape parameter of RBFs as shown in \cite{lee2003local}. These localized meshless methods have been successfully applied to a large variety of problems to solve partial differential equations including the Richards equation \cite{stevens2009order,stevens2011scalable,ben2018radial,boujoudar2021localized,boujoudar2022modelling}. In this study, LRBF methods are applied to solve the resulting Kirchhoff-transformed Richards equation for modeling infiltration through soils.

The paper is organized as follows. In Section \ref{sec:2}, we introduce the developed numerical approach and the Kirchhoff-transformed Richards equation. The proposed numerical model based on the LRBF method is described in Section \ref{sec:3}. In Section \ref{sec:4}, numerical simulations are conducted to validate the developed approach for modeling infiltration through soils. Finally, some concluding remarks are provided in Section \ref{sec:5}.

\section{Material and models}
\label{sec:2}

\subsection{Richards' model} 
\label{sec:2.1}
We consider the traditional Richards equation describing infiltration through soils \cite{richards1931capillary}:
  \begin{equation}\label{E2.1}
  \dfrac{\partial\theta(h)}{\partial t}-\nabla.\left(K_\text{s}(\boldsymbol{x})k_\text{r}(h)\nabla (h+z)\right)= s(\boldsymbol{x},t), \text{ $\boldsymbol{x}\in\Omega$},
  \end{equation}
where $\theta$ $ [L^{3}/L^{3}]$ is the water content, $h$ $[L]$ is the pressure head, $K_\text{s}$ $[L/T]$ is the saturated hydraulic conductivity which depends on the medium's spatial heterogeneity, $k_\text{r}$ $[-]$ is the water relative permeability, $s( \boldsymbol{x},t)$ is a source or sink term which may include evaporation and plant-root extraction, $\Omega $ denotes an open subset of $\mathbb{R}^{3}$, $\boldsymbol{x}$ $[L]$ is the spatial coordinate and $z$ $[L]$ is the upward vertical coordinate.\\ We consider the water-saturation $S=(\theta-\theta_{r})/(\theta_{s}-\theta_{r})$ $[-]$ to write Richards' equation without source/sink term: 
\begin{equation} \label{E2.3}
\phi \dfrac{\partial S}{\partial t}-\nabla.(K_\text{s}(\boldsymbol{x})k_\text{r}(S)\nabla h)-\dfrac{\partial (K_\text{s}k_\text{r})}{\partial z}= 0, \text{ $\boldsymbol{x}\in\Omega$},
\end{equation}  
where $\theta_{s}$ $ [L^{3}/L^{3}]$ is the saturated water content, $\theta_{r}$ $ [L^{3}/L^{3}]$ is the residual water content
and the parameter $\phi$ $[-]$ is given by $\phi=(\theta_{s}-\theta_{r})$.

Equation (\ref{E2.3}) is highly non-linear due to the nonlinear dependencies of the capillary pressure and relative permeability functions on saturation. Empirical constitutive relationships have been developed for these functions using experiments  \cite{gardner1958some,brooks1964hydraulic,van1980closed} and are used in previous studies \cite{celia1990general, cueto2008nonlocal, beljadid2020continuum, boujoudar2021localized, keita2021implicit}. 
\subsection{Capillary pressure and relative permeability functions}
\label{sec:2.2}
 
Here, we used the Brooks and Corey's model~\cite{brooks1964hydraulic} for the capillary  pressure function. The saturation is given by:
 \begin{equation} \label{E2.A2}
S(h)=
\begin{cases}
\left( \dfrac{h}{h_{\text{d}}}\right)^{-\lambda},  &\text{if}~  h \leqslant h_{\text{d}}, \\
1, &\text{if}~ h > h_{\text{d}},
\end{cases}
\end{equation}
where $\lambda$ $[-]$ is the Brooks-Corey parameter,  $h_{\text{d}}=-h_{\text{cap}}$ $[L]$ and the characteristic capillary rise $h_{\text{cap}}$ is determined by the Leverett scaling formula \cite{leverett1941capillary}.
The capillary pressure function can be expressed as follows:
 \begin{equation}\label{EJ.6}
 h(S)=h_{\text{d}}J(S),
 \end{equation}
 where the Leverett $J$-function is given by $J(S)=S^{-1/\lambda}$. We propose to use the power-law relation in saturation for the relative permeability $k_{r}(S)=S^{\beta}$ which can be expressed using capillary pressure as follows:
\begin{equation} \label{E2.7}
k_{r}(h)=\begin{cases}
\left( \dfrac{h}{h_{\text{d}}}\right)^{-\lambda \beta},  &\text{if}~  h  \leqslant h_{\text{d}} \\
1, &\text{if}~ h > h_{\text{d}},
\end{cases}
\end{equation}
 where $\beta > 1$.\\
 
Let's introduce a reference constant $\bar{h}~[L]$ for the capillary pressure in Equation (\ref{E2.3}) in order to use a dimensionless form of the capillary pressure function. For $ h\leqslant h_{\text{d}}$, we obtain the following expression where we set $\omega= \bar{h}/h_{\text{d}}$:
\begin{equation}\label{E2.10}
K_{\text{s}}k_{\text{r}}\nabla h = K_{\text{s}}\omega^{-\lambda \beta} \left( \dfrac{h}{\bar{h}}\right)^{-\lambda \beta}\nabla h, 
\end{equation}
and $\bar{h}$ can be taken as:
\begin{equation}
\bar{h}=\dfrac{1}{V}\iiint_{V} h_{\text{d}}(x,y,z)dx dy dz,
\end{equation}
where $V$ is the volume of the whole domain.
Since $\omega$ depends only on space ($\omega=\omega(\boldsymbol{x})$), we obtain:
\begin{equation}
 \dfrac{\partial S}{\partial t} =\omega^{-\lambda }\frac{\partial  }{\partial t}\left[ \left( \dfrac{h}{\bar{h}}\right)^{-\lambda}\right] ,\label{E2.13}
\end{equation}
and
\begin{equation}\label{E2.161}
 \dfrac{\partial (K_\text{s}k_\text{r})}{\partial z} =\dfrac{\partial }{\partial z}\left[ K_\text{s} \omega^{-\lambda \beta} \left( \dfrac{h}{\bar{h}}\right)^{-\lambda \beta}\right].
\end{equation}
We substitute Equations (\ref{E2.10}), (\ref{E2.13}) and (\ref{E2.161}) into Equation (\ref{E2.3}), we obtain the new form of the Richards equation associated with the Brooks and Corey model for $ h\leqslant h_{\text{d}}$:
\begin{equation}
\phi\omega^{-\lambda }\frac{\partial  }{\partial t}\left[ \left( \dfrac{h}{\bar{h}}\right)^{-\lambda}\right]- \nabla.\left[ K_{\text{s}}\omega^{-\lambda \beta} \left( \dfrac{h}{\bar{h}}\right)^{-\lambda \beta}\nabla h\right]-\dfrac{\partial }{\partial z}\left[ K_\text{s} \omega^{-\lambda \beta} \left( \dfrac{h}{\bar{h}}\right)^{-\lambda \beta}\right] =0,\label{E2.17}
\end{equation}
where we separate the terms which are dependent on the capillary pressure function and those which depend only on space due to soil heterogeneity.\\
For $ h > h_{\text{d}}$, we obtain:
\begin{equation}\label{E27}
-\nabla.(K_{\text{s}}\nabla h)-\dfrac{\partial K_\text{s}}{\partial z}=0.
\end{equation}
\subsection{Kirchhoff transformation}
\label{sec:2.3}
In this study, we propose to use the transformation of Kirchhoff which allows us to reduce the nonlinearity of the model equation. The Kirchhoff integral transformation is defined as:
\begin{equation}\label{KT}
\varphi(h) =\bar{h}\int_{+\infty}^{ h/\bar{h}}\sigma^{-\lambda\beta}d\sigma,
\end{equation}
where we used the variable $\sigma=h/\bar{h}$. Based on Equation (\ref{E2.7}), the transformation (\ref{KT}) can be rewritten as follows:
\begin{equation}\label{E2.14} 
\varphi(h)=
\begin{cases}
\dfrac{\bar{h}}{(1-\lambda\beta)}\left( \dfrac{h}{\bar{h}}\right)^{(1-\lambda \beta)},  &\text{if}~  h \leqslant h_{\text{d}}, \\
\\
\dfrac{\bar{h}}{(1-\lambda\beta)}\left( \dfrac{h_{\text{d}}}{\bar{h}}\right)^{(1-\lambda \beta)}+\left( \dfrac{h_{\text{d}}}{\bar{h}}\right)^{-\lambda\beta}\left( h-h_{\text{d}}\right), &\text{if}~ h > h_{\text{d}},
\end{cases}
\end{equation}
where, we assume that $\lambda\beta>1$. From Equation (\ref{E2.14}), we can determine the pressure head $h$ as a function of $\varphi$: 
\begin{equation}\label{E2.16} 
h=
\begin{cases}
\bar{h}\left( \dfrac{(1-\lambda\beta)}{\bar{h} }\varphi\right) ^{1/(1-\lambda\beta)},  &\text{if}~  \varphi \leqslant   \dfrac{\bar{h}}{(1-\lambda\beta)}\left( \dfrac{h_{\text{d}}}{\bar{h}}\right)^{(1-\lambda \beta)}, \\
\\
\left( \dfrac{h_{\text{d}}}{\bar{h}}\right)^{\lambda \beta}\varphi+h_{\text{d}}-\dfrac{h_{\text{d}}}{1-\lambda\beta}, &\text{if}~ \varphi > \dfrac{\bar{h}}{(1-\lambda\beta)}\left( \dfrac{h_{\text{d}}}{\bar{h}}\right)^{(1-\lambda \beta)}.
\end{cases}
\end{equation}

In our approach, in the transformation of the model equation, the variation of the intrinsic permeability is assumed dominant in the effect of heterogeneity variability \cite{leverett1941capillary, keita2021implicit}. We have for $h \leqslant h_{\text{d}}$:
\begin{equation}
\nabla\varphi(h) =\left( \dfrac{h}{\bar{h}}\right)^{-\lambda\beta} \nabla h,
\end{equation}
and
\begin{equation}
\nabla.\left[ K_{\text{s}}\omega^{-\lambda \beta} \left( \dfrac{h}{\bar{h}}\right)^{-\lambda \beta}\nabla h\right]=\nabla.\left(  K_{\text{s}}\omega^{-\lambda \beta} \nabla\varphi\right),
\end{equation}\label{E2.21}
and the third term of Equation (\ref{E2.17}) becomes:
\begin{align}
\dfrac{\partial }{\partial z}\left[ K_\text{s} \omega^{-\lambda \beta} \left( \dfrac{h}{\bar{h}}\right)^{-\lambda \beta}\right]&=\dfrac{\partial }{\partial z}\left[ K_\text{s} \omega^{-\lambda \beta} \dfrac{(1-\lambda\beta)}{\bar{h}}\left( \dfrac{h}{\bar{h}}\right)^{-1} \varphi\right].
\end{align}
The time derivative in the first term of Equation (\ref{E2.17}) can be rewritten under the same assumption $h \leqslant h_{\text{d}}$ by:
\begin{equation}
\frac{\partial  }{\partial t}\left[ \left( \dfrac{h}{\bar{h}}\right)^{-\lambda}\right]=\frac{-\lambda }{\bar{h}} \left( \dfrac{h}{\bar{h}}\right)^{-\lambda-1} \dfrac{\partial h }{\partial t},
\end{equation}
and
\begin{equation}
\frac{\partial \varphi }{\partial t}= \left( \dfrac{h}{\bar{h}}\right)^{-\lambda \beta}\frac{\partial h }{\partial t},
\end{equation}
which implies:
\begin{equation}
\frac{\partial  }{\partial t}\left[ \left( \dfrac{h}{\bar{h}}\right)^{-\lambda}\right]=\frac{-\lambda}{\bar{h}}\left( \dfrac{h}{\bar{h}}\right)^{\lambda\beta-\lambda-1}\frac{\partial \varphi }{\partial t}.
\end{equation}

Similarly, for $ h > h_{\text{d}}$, Equation (\ref{E27}) can be written in terms of $\varphi$:
\begin{equation}\label{E2.283}
-\nabla. (K_{\text{s}}\omega^{-\lambda \beta} \nabla\varphi)-\dfrac{\partial K_{\text{s}} }{\partial z}=0.
\end{equation}
For simplicity, we will use the following parameters:
\begin{equation}
\chi=K_{\text{s}}\omega^{-\lambda \beta},
\end{equation}
\begin{equation}
E=
\begin{cases}
\phi\dfrac{-\lambda}{\bar{h}}\omega^{-\lambda }\left( \dfrac{h}{\bar{h}}\right)^{\lambda\beta-\lambda-1}, & \text{if}~h \leqslant h_{\text{d}},\\
0, & \text{if}~h > h_{\text{d}},
\end{cases}
\end{equation}
\begin{equation}
F=
\begin{cases}
\dfrac{(1-\lambda\beta)}{\bar{h}}\left( \dfrac{h}{\bar{h}}\right)^{-1}, & \text{if}~h \leqslant h_{\text{d}},\\
0, & \text{if}~h > h_{\text{d}},
\end{cases}
\end{equation}
\begin{equation}
G=
\begin{cases}
0, & \text{if}~h \leqslant h_{\text{d}},\\
K_{\text{s}}, & \text{if}~h > h_{\text{d}}.
\end{cases}
\end{equation}

Finally, we obtain the new form of the model equation using the Kirchhoff $\varphi$:
\begin{equation}\label{E2.281}
E\frac{\partial \varphi }{\partial t}-\nabla. (\chi \nabla\varphi)-\dfrac{\partial }{\partial z}\left(\chi F \varphi\right)-\dfrac{\partial G }{\partial z}=0,
\end{equation}
where the non-linearity of the original Richards model has been reduced since only the terms $E$ and $F$ are nonlinear and $\chi$ depends only on spatial coordinates $\boldsymbol{x}$ due to the heterogeneity of soils.
\section{Numerical model} \label{sec:3}
\subsection{ Approximation methods}
In this section, we describe the approximation methods used to solve Equation (\ref{E2.281}). Let $\Delta t>0$ a time step and $t^{p}=p\Delta t$ with $p\geq0$ denotes the time level. Temporal discretization of Equation (\ref{E2.281}) using the backward Euler method may be written as:
\begin{equation}\label{E2.82}
E^{p+1}\dfrac{\varphi^{p+1}-\varphi^{p}}{\Delta t}-\nabla.(\chi \nabla\varphi^{p+1})-\dfrac{\partial \left(\chi F^{p+1} \varphi^{p+1} \right)}{\partial z}-\dfrac{\partial G }{\partial z}=0,
\end{equation}
where $\varphi^{p+1}$ denotes the approximate solution of $\varphi$ at $t^{p+1}$, $E^{p+1}$ and $F^{p+1}$ are the estimated values of $E$ and $F$ computed using $h^{p+1}$ which is obtained by substituting $\varphi^{p+1}$ in Equation (\ref{E2.16}). 

By applying the Picard iteration scheme to Equation (\ref{E2.82}), we obtain:
\begin{equation}\label{E2.282}
E^{p+1,m}\dfrac{\varphi^{p+1,m+1}-\varphi^{p}}{\Delta t}-\nabla.(\chi \nabla\varphi^{p+1,m+1})-\dfrac{\partial }{\partial z}\left(\chi  F^{p+1,m} \varphi^{p+1,m+1} \right) -\dfrac{\partial G }{\partial z}=0,
\end{equation}
with $m$ identifies iteration level. The solution is assumed to be known both at time level $p$  and
at iteration level $m$.
Let $\left\lbrace \boldsymbol{x_{i}}=(x_{i},y_{i},z_{i}) \right\rbrace _{i=1}^{N_{i}}\subset\Omega$ be $N_{i}$ uniform distinct points  and $\left\lbrace \boldsymbol{x_{i}}\right\rbrace _{i=N_{i}+1}^{N}\subset\partial\Omega$ be $N_{b}$ distinct nodes, where $N_{i}$ denotes the number of interior points and $N_{b}$ denotes the number of points on the boundary $(N=N_{i}+N_{b})$. 

For each point $\left( \boldsymbol{x_{i}} \right)_{i=1}^{N_{i}}$, we discretize $\nabla.(\chi \nabla\varphi^{p+1,m+1})$ as follows:
\begin{equation}
\nabla.(\chi \nabla\varphi^{p+1,m+1})=\dfrac{\partial }{\partial x}\left( \chi\dfrac{\partial\varphi^{p+1,m+1}}{\partial x}   \right)+\dfrac{\partial }{\partial y}\left(  \chi \dfrac{\partial\varphi^{p+1,m+1}}{\partial y}  \right)+\dfrac{\partial }{\partial z}\left(  \chi \dfrac{\partial\varphi^{p+1,m+1}}{\partial z}  \right).
\end{equation}
For simplicity, we take the following expressions:
\begin{equation}
 \mathcal{L}_{d}^{m}\varphi_{i}=\dfrac{\partial }{\partial x^{(d)}}\left( \chi_{i}\dfrac{\partial\varphi_{i}}{\partial x^{(d)}} \right),   
\end{equation}
\begin{equation}
 \mathcal{L}_{4}^{m}\varphi_{i}=\dfrac{\partial }{\partial z}\left(\chi_{i} F_{i}  \varphi_{i} \right), 
\end{equation}
where $d=\left\lbrace 1,2,3\right\rbrace $ and $ (x^{(1)},x^{(2)},x^{(3)})=(x,y,z) $. By extending the $1D$ spatial approximation used in \cite{celia1990general} to $3D$ case, we obtain:
\begin{equation}\label{eqA}
\mathcal{L}_{d}^{m}\varphi_{i}=\dfrac{1}{(\Delta x^{(d)})^{2}}\left(\chi^{(d)}_{i+1/2}(\varphi^{(d)}_{iR}-\varphi^{(d)}_{i})-\chi^{(d)}_{i-1/2}(\varphi^{(d)}_{i}-\varphi^{(d)}_{iL})\right),
\end{equation}
 where the expressions $\chi^{(d)}_{i+1/2}$ and $\chi^{(d)}_{i-1/2}$ are given by:
\begin{equation}
\begin{cases}
\chi^{(d)}_{i+1/2}=\dfrac{1}{2}(\chi^{(d)}_{i}+\chi^{(d)}_{iR}),\\
\chi^{(d)}_{i-1/2}=\dfrac{1}{2}(\chi^{(d)}_{i}+\chi^{(d)}_{iL}).
\end{cases}
\end{equation}
Along the $x^{(d)}$-axis, $\varphi^{(d)}_{iR}$ and $\chi^{(d)}_{iR}$ are the corresponding values at the right of $\varphi_{i}$ and $\chi_{i}$ respectively. Similarly, $\varphi^{(d)}_{iL}$ and $\chi^{(d)}_{iL}$ are the corresponding values at the left of $\varphi_{i}$ and $\chi_{i}$ respectively. 
We then obtain:
\begin{equation}
\nabla.(\chi \nabla\varphi^{p+1,m+1})=\mathcal{L}_{1}^{m}\varphi^{p+1,m+1}+\mathcal{L}_{2}^{m}\varphi^{p+1,m+1}+\mathcal{L}_{3}^{m}\varphi^{p+1,m+1}.
\end{equation}
On the other hand, we approximate the following spatial operators as follows \cite{celia1990general}:
\begin{equation}
\mathcal{L}_{4}^{m}\varphi_{i}=\dfrac{1}{(\Delta z)}\left(\chi^{(3)}_{i+1/2}F^{(3)}_{i+1/2}\varphi^{(3)}_{i+1/2}-\chi^{(3)}_{i-1/2}F^{(3)}_{i-1/2} \varphi^{(3)}_{i-1/2} \right),
\end{equation}
\begin{equation}
\dfrac{\partial G_{i} }{\partial z}=\dfrac{1}{( \Delta z)}\left(G^{(3)}_{i+1/2}-G^{(3)}_{i-1/2} \right).
\end{equation}
We then obtain:
\begin{equation}
\dfrac{\partial }{\partial z}\left(\chi_{i} F^{p+1,m}_{i}  \varphi^{p+1,m+1}_{i} \right)=\mathcal{L}_{4}^{m}\varphi^{p+1,m+1}_{i}.
\end{equation}

For simplicity, we use the following expressions:
\begin{equation}
\mathcal{L}^{m}\varphi^{p+1,m+1}_{i}=\dfrac{E^{p+1,m}}{\Delta t}\varphi^{p+1,m+1}_{i}-\left( \mathcal{L}_{1}^{m}.+\mathcal{L}_{2}^{m}.+\mathcal{L}_{3}^{m}.\right)\varphi^{p+1,m+1}_{i}-\mathcal{L}_{4}^{m}\varphi^{p+1,m+1}_{i},
\end{equation}
\begin{equation}
f_{i}^{p+1,m}=\dfrac{E^{p+1,m}}{\Delta t}\varphi_{i}^{p}+\dfrac{1}{\Delta z}(G^{(3)}_{i+1/2}-G^{(3)}_{i-1/2}).
\end{equation}
The operator $\mathcal{L}^{m}$ is linear for each iteration level $m$.
In addition to initial and boundary conditions, Equation (\ref{E2.282}) may be rewritten as follows: 
\begin{equation}\label{E3.5}
\begin{cases}
\mathcal{L}^{m}\varphi^{p+1,m+1}(\boldsymbol{x})=f^{p+1,m}(\boldsymbol{x}), & \text{ $\boldsymbol{x}\in\Omega$},\\
\mathcal{B}\varphi^{p+1,m+1}(\boldsymbol{x})=\varphi_{\Gamma}(\boldsymbol{x}), & \text{ $\boldsymbol{x}\in\partial\Omega$},\\
\varphi^{0,m+1}(\boldsymbol{x})=\varphi_{0}(\boldsymbol{x}), & \text{ $\boldsymbol{x}\in\Omega$.}
\end{cases}
\end{equation}
$\varphi_{0}$ and $\varphi_{\Gamma}$ are functions associated with the initial and boundary conditions. For each time level $p$, the linear system (\ref{E3.5}) is solved at each iteration level of Picard $m$ until the following inequality is satisfied at all collocation points:
\begin{equation}
\delta^{m}=\lvert \varphi^{p+1,m+1}-\varphi^{p+1,m}\rvert\leq Tol,
\end{equation}
where $Tol$ is the error tolerance. 
\subsection{Approach using local radial basis functions}\label{s23}
In this study, we use the local RBF meshfree method \cite{li2013localized} which has advantages in terms of operational memory and calculation procedures where only inversions of small size matrices are required.\\For any point $\boldsymbol{x_{s}}\in\bar{\Omega}$, the \textit{k}-d tree algorithm is used \cite{bentley1975multidimensional} to create a localized influence domain $\Omega^{[s]}=\left\lbrace \boldsymbol{x_{k}}^{[s]} \right\rbrace _{k=1}^{n_{s}}\subset\bar{\Omega}$. It contains $n_{s}$ nearest nodal points from $\boldsymbol{x_{s}}$.
\\
In the local RBF approach, the transformed Kirchhoff variable $\varphi_{[s]}^{p+1,m+1}$ is approximated in each localized influence domain $\Omega^{[s]} $ as follows: 
\begin{equation}\label{E3.15}
\varphi_{[s]}^{p+1,m+1}(\boldsymbol{x_{s}})=\sum_{i=1}^{n_{s}}\alpha^{p+1,m+1}_{i}\psi(\Vert \boldsymbol{x_{s}}-\boldsymbol{x_{i}}^{[s]} \Vert ),
\end{equation}
where $ \lbrace \alpha_{i}^{p+1,m+1} \rbrace_{i=1}^{n_{s}} $ are constants to be determined and $ \psi $ is a RBF for which, in our approach, we chose the exponential function given by $ \psi(r)=\exp({-(c r)^{2}})$, where $ r=\Vert \boldsymbol{x_{s}} -\boldsymbol{x_{i}}^{[s]} \Vert$ denotes the distance between $\boldsymbol{x_{s}}$ and $\boldsymbol{x_{i}}^{[s]}$ and $c>0 $ is the shape parameter.

According to Equation~(\ref{E3.15}), we obtain:
\begin{equation}\label{E3.16d}
\varphi_{[s]}^{p+1,m+1}=\psi^{[s]} \alpha^{p+1,m+1}_{[s]},
\end{equation}
where we used the matrix $\psi^{[s]}=\left[ \psi(\Vert \boldsymbol{x_{i}}^{[s]}-\boldsymbol{x_{j}}^{[s]} \Vert) \right]_{1\leqslant i,j \leqslant n_{s}} $ and the vectors: \\
$\varphi_{[s]}^{p+1,m+1}=\left[ \varphi_{[s]}^{p+1,m+1}(\boldsymbol{x_{1}}^{[s]}),\varphi_{[s]}^{p+1,m+1}(\boldsymbol{x_{2}}^{[s]}),...,\varphi_{[s]}^{p+1,m+1}(\boldsymbol{x_{n_{s}}}^{[s]}) \right]^{T}$, \\$\alpha^{p+1,m+1}_{[s]}=\left[ \alpha^{p+1,m+1}_{1},\alpha^{p+1,m+1}_{2},...,\alpha^{p+1,m+1}_{n_{s}} \right]^{T} $. 

From Equation (\ref{E3.16d}), we obtain:
\begin{equation}
\alpha^{p+1,m+1}_{[s]}=(\psi^{[s]})^{-1}\varphi_{[s]}^{p+1,m+1}.
\end{equation}
Applying the linear operator $\mathcal{L}^{m}$ to Equation~(\ref{E3.15}) at each $\boldsymbol{x_{s}} \in \Omega^{[s]}$, we have:

\begin{equation} \label{E3.18}
\begin{split}
\mathcal{L}^{m}\varphi_{[s]}^{p+1,m+1}(\boldsymbol{x_{s}}) & = \sum_{i=1}^{n_{s}}\alpha^{p+1,m+1}_{i}\mathcal{L}^{m}\psi(\Vert \boldsymbol{x_{s}}-\boldsymbol{x_{i}}^{[s]} \Vert )=\sum_{i=1}^{n_{s}}\alpha^{p+1,m+1}_{i}\Psi^{m}(\Vert \boldsymbol{x_{s}}-\boldsymbol{x_{i}}^{[s]} \Vert ) \\    
& = \Gamma_{[s]}^{m}\alpha^{p+1,m+1}_{[s]}=\Gamma_{[s]}^{m}(\psi^{[s]})^{-1}\varphi_{[s]}^{p+1,m+1}=\Upsilon_{[s]}^{m}\varphi_{[s]}^{p+1,m+1},
\end{split}
\end{equation}
where $ \Psi^{m}=\mathcal{L}^{m}\psi$, $\Gamma_{[s]}^{m}=\left[ \Psi(\Vert \boldsymbol{x_{s}}-\boldsymbol{x_{1}}^{[s]} \Vert),...,\Psi(\Vert \boldsymbol{x_{s}}-\boldsymbol{x_{n_{s}}}^{[s]} \Vert) \right]$ and $ \Upsilon_{[s]}^{m}=\Gamma_{[s]}^{m}(\psi^{[s]})^{-1} $.
\\
In order to reformulate Equation (\ref{E3.18}) in terms of the global vector $\varphi^{p+1,m+1}$ instead of $\varphi_{[s]}^{p+1,m+1}$, $\Upsilon^{m}$ is considered as the expansion of $\Upsilon_{[s]}^{m}$ by inserting zeros in the proper position. It follows that:
\begin{equation}\label{E3.71}
\mathcal{L}^{m}\varphi_{[s]}^{p+1,m+1}(\boldsymbol{x_{s}})= \Upsilon^{m} \varphi^{p+1,m+1}, 
\end{equation}
where $ \varphi^{p+1,m+1}=\left[ \varphi^{p+1,m+1}(\boldsymbol{x_{1}}),\varphi^{p+1,m+1}(\boldsymbol{x_{2}}),...,\varphi^{p+1,m+1}(\boldsymbol{x_{N}}) \right]^{T} $.\\
Similarly, for $ \boldsymbol{x_{s}}\in\partial\Omega$, we apply the linear operator $\mathcal{B}$:
\begin{equation} \label{E3.17}
\begin{split}
\mathcal{B}\varphi_{[s]}^{p+1,m+1}(\boldsymbol{x_{s}})&=\sum_{i=1}^{n_{s}}\alpha^{p+1,m+1}_{i}\mathcal{B}\psi(\Vert \boldsymbol{x_{s}}-\boldsymbol{x_{i}}^{[s]} \Vert  )=(\mathcal{B}\psi^{[s]})\alpha^{p+1,m+1}_{[s]}\\
&=(\mathcal{B}\psi^{[s]})(\psi^{[s]})^{-1}\varphi_{[s]}^{p+1,m+1}=\mathbb{\upsilon}^{[s]}\varphi_{[s]}^{p+1,m+1}=\mathbb{\upsilon}\varphi^{p+1,m+1},
\end{split}
\end{equation}
where $\mathbb{\upsilon}^{[s]}=(\mathcal{B}\psi^{[s]})(\psi^{[s]})^{-1} $ and  $ \mathbb{\upsilon}$ is the global expansion of $ \mathbb{\upsilon}^{[s]}$ by adding zeros in the proper location.

From Equations (\ref{E3.71}) and (\ref{E3.17}), we
get the system below:
\begin{equation}
\begin{split}
\mathcal{L}^{m}\varphi^{p+1,m+1}(\boldsymbol{x_{s}})=\Upsilon^{m}(\boldsymbol{x_{s}})\varphi^{p+1,m+1}=f^{p+1,m}(\boldsymbol{x_{s}}),\\
\mathcal{B}\varphi^{p+1,m+1}(\boldsymbol{x_{s}})=\mathbb{\upsilon}(\boldsymbol{x_{s}})\varphi^{p+1,m+1}=\varphi_{\Gamma}(\boldsymbol{x_{s}}).
\end{split}
\end{equation}
We obtain the following sparse linear system: 
\begin{equation}\label{3.20}
\left(
\begin{array}{c}
\Upsilon^{m}(\boldsymbol{x_{1}}) \\
\Upsilon^{m}(\boldsymbol{x_{2}}) \\
 . \\
.  \\
\Upsilon^{m}(\boldsymbol{x_{N_{i}}})\\
\mathbb{\upsilon}(\boldsymbol{x_{N_{i}+1}}) \\
 . \\
.  \\
\mathbb{\upsilon}(\boldsymbol{x_{N}})
\end{array}
\right)
\left(
\begin{array}{c}
\varphi^{p+1,m+1}(\boldsymbol{x_{1}}) \\
\varphi^{p+1,m+1}(\boldsymbol{x_{2}})\\
.\\
.\\
\varphi^{p+1,m+1}(\boldsymbol{x_{N_{i}}})\\
\varphi^{p+1,m+1}(\boldsymbol{x_{N_{i}+1}})\\
.\\
.\\
\varphi^{p+1,m+1}(\boldsymbol{x_{N}})
\end{array}
\right)
=
\left(
\begin{array}{c}
f^{p+1,m}(\boldsymbol{x_{1}})\\
f^{p+1,m}(\boldsymbol{x_{2}})\\
.\\
.\\
f^{p+1,m}(\boldsymbol{x_{N_{i}}})\\
\varphi_{\Gamma}(\boldsymbol{x_{N_{i}+1}})\\
. \\
.\\
\varphi_{\Gamma}(\boldsymbol{x_{N}})
\end{array}
\right).
\end{equation}
The LRBF approach leads to system of sparse equations (\ref{3.20}) which allows as to reduce the size of the dense matrices  and avoid ill-conditioned problems arising from the global approach \cite{yao2012assessment,young2016localized,hamaidi2016space,li2022efficient}. The approximate solutions $\varphi^{p+1,m+1}=\left\lbrace\varphi^{p+1,m+1}(\boldsymbol{x}_{i} )\right\rbrace_{i=1}^{N}$ can be obtained by solving the system (\ref{3.20}).

 \subsection{Initial and boundary conditions}\label{IBC}
At $t=0$, we assume that $h(\boldsymbol{x},0)=h_{0}(\boldsymbol{x})$ which implies in terms of the Kirchhoff variable that $\varphi(\boldsymbol{x},0)=\varphi_{0}(\boldsymbol{x})$, where $$\varphi_{0}(\boldsymbol{x})=\dfrac{\bar{h}}{(1-\lambda\beta)}\left( \dfrac{h_{0}}{\bar{h}}\right) ^{(1-\lambda \beta)}.$$ The boundary conditions are expressed in terms of the Kirchhoff variable. For Dirichlet conditions, we consider $h=g_{b}$ for $z = 0$ and $z = L$ which implies that: $$\varphi_{b}= \dfrac{\bar{h}}{(1-\lambda \beta)} \left( \dfrac{g_{b}}{\bar{h}}\right)^{(1-\lambda \beta)},$$ where $g_{b}$ is given by: $$ g_{b}=\begin{cases} h_{0} , &  z=0,\\ 0, &  z=L. \end{cases}$$ For Neumann conditions, we consider $-K\dfrac{\partial h} {\partial x}=0$ for $x=0$ and $x=l_{1}$ which implies that $-\chi\dfrac{\partial \varphi} {\partial x}=0$. In the same way, we assume that $-\chi\dfrac{\partial \varphi} {\partial y}=0$ at the lateral boundaries $y = 0$ and $y =l_{2}$. Therefore, the linear operator corresponding to the boundary conditions is given by:
\begin{equation}\label{E328}
\mathcal{B}\varphi=\begin{cases}
\varphi_{b} , &  z=\left\{0,L\right\},\\
-\chi\dfrac{\partial \varphi} {\partial x}, & x=\left\{0,l_{1}\right\}, \\
-\chi\dfrac{\partial \varphi} {\partial y}, & y=\left\{0,l_{2}\right\}.
\end{cases}
\end{equation}.

\section{\textbf{Numerical experiments}}\label{sec:4}
To validate the proposed approach for modeling unsaturated flow through heterogeneous soils, we present numerical solutions of Richards equation in one-, two- and three-dimensional systems. The computational domain $\Omega=[0,l_{1}]\times[0,l_{2}]\times[0,L]$ is used to perform $3D$ numerical simulations and we consider the domains $\Omega=[0,l_{1}]\times[0,L]$ and $\Omega=[0,L]$, respectively, for $2D$ and $1D$ numerical simulations.\\
In Section (\ref{ss41}), numerical tests are performed using the developed model to simulate flow in unsaturated homogeneous soils. Section (\ref{ss42}) presents numerical simulations of infiltration through heterogeneous soils.
To investigate the accuracy of the developed technique, the results of numerical tests are used to compute the $RMSE$ and $L_{er}^{1}$ errors based on the following formulas:
\begin{equation}
RMSE=\sqrt{\dfrac{1}{N}\sum_{i=1}^{N}\vert \theta(x_{i})-\theta_{ref}(x_{i}) \vert^{2}},
\end{equation}
\begin{equation}
L_{er}^{1}=\dfrac{\sum_{i=1}^{N}( \theta(x_{i})-\theta_{ref}(x_{i}))^{2} }{\sum_{i=1}^{N} \theta_{ref}(x_{i})^{2} },
\end{equation}
where $\theta(x_{i})$ represents the approximate solution for the water content and $\theta_{ref}(x_{i})$ represents a reference solution and $N$ is the number of collocation points.

\subsection{Infiltration in homogeneous soils}\label{ss41}
\begin{figure}[ht]
\centering
\begin{tabular}{cc}
Clay & Clay loam \\
\includegraphics[width=6.2cm,height=6cm,angle=0]{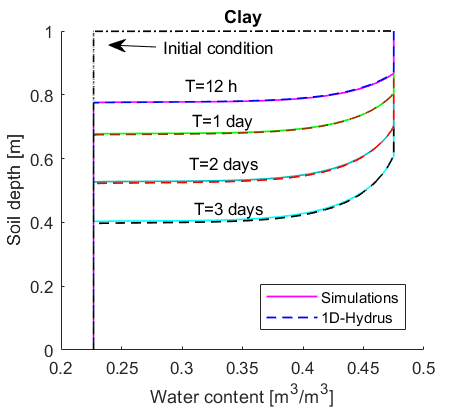} & 
\includegraphics[width=6.2cm,height=6cm,angle=0]{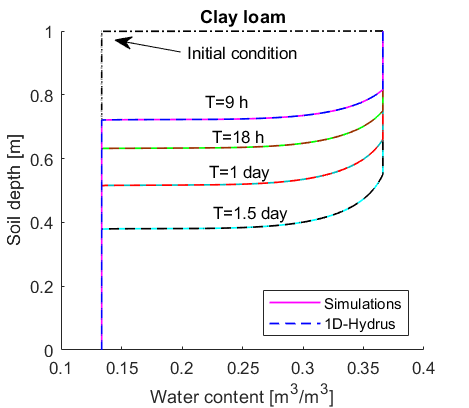} \\
Sand & Silty clay \\
\includegraphics[width=6.2cm,height=6cm,angle=0]{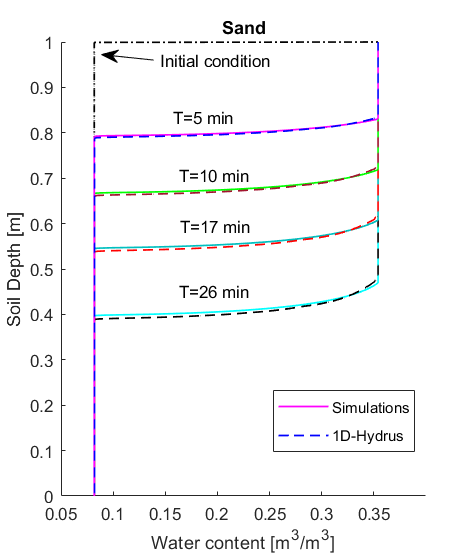} &  
\includegraphics[width=6.2cm,height=6cm,angle=0]{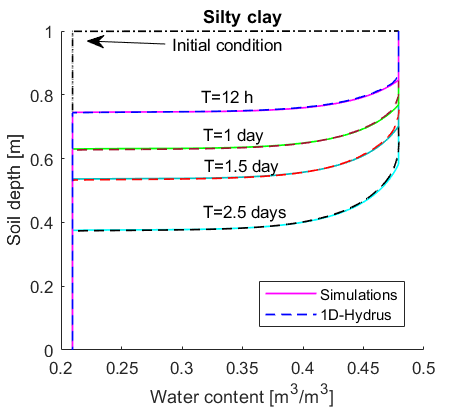} \\
\end{tabular}
\caption{The water content profiles of the numerical and the reference solutions.}\label{f1}
\end{figure}

In this section, we perform numerical simulations using the numerical method to simulate unsaturated flow in homogeneous soils. We consider different soil samples with $L=1~m$. The parameters of these soils are shown in Table \ref{t1}.
\begin{table}[h]
\begin{center}
\caption{Parameters of soils. }\label{t1}
 \begin{tabular}{|c|c|c|c|c|c|c|c|c|}
 \hline
Soil& Type  &$\theta_r$ & $\theta_s$   &   $\theta_0$ &$K_{s}$ & $h_{\text{d}}$  & $\lambda$ & $\beta$ \\
 --& -- & $(m^3/m^3)$  &   $(m^3/m^3)$ & $(m^3/m^3)$ & $(m/day)$ & $(m)$ &-- & --\\ 
 \hline
 1 &  Clay  & $0.09$ & $0.475$ & $0.226$ & $0.0144$ & $-0.3731$  & $0.131$ &$18.2672$ \\
 \hline
 2 & Clay loam & $0.075$ & $0.366$ & $0.130$ & $0.040$ & $-0.2590$ &  $0.194$ &$13.3093$\\
 \hline
 3 & Sand& $0.04$ & $0.354$ & $0.0819$ & $5.04$ & $-0.01471$  & $1.051$ &$ 4.9029$ \\
 \hline
 4 & Silty clay & $0.056$ & $0.479$ & $0.212$ & $0.0216$ & $-0.3425$ & $0.127$ & $ 18.7480$ \\
 \hline
  \end{tabular}
\end{center}
 \end{table} 
 The parameters $\theta_0$ and $h_0$ are the initial water content and pressure head respectively. In this numerical test, we used $N_{z} = 1001$ uniform nodes with $n_{s} =3$ number of neighboring points and a time step $\Delta t=0.0001$.
The water content profiles are shown in Figure \ref{f1} for the considered types of soils.
Table \ref{3DRMSE1} illustrates the $RMSE$, $L^{1}_{er}$ errors between the numerical solutions and the reference solutions which are obtained using $1D$-Hydrus \cite{simunek2005hydrus}. We obtain accurate results and the predictions are in good agreement with the reference solutions simulated using $1D$-Hydrus.

\begin{table}[ht]
\begin{center}
\caption{The computed errors.}\label{3DRMSE1}
 \begin{tabular}{|c|c|c|c|}
 \hline
 Soils & $T$ & $RMSE$ & $L^{1}_{er}$ \\
 \hline
 \multirow{2}{1.5cm}{1} & $12$ h& $1.2
 \times 10^{-3}$ & $8.31\times 10^{-4}$  \\
 \cline{2-4}
 & $3~$days  & $7.7
 \times 10^{-3}$   & $3 \times 10^{-3}$  \\
 \cline{2-4}
 \hline
 \multirow{2}{1.5cm}{2} & $9$ h  & $6.4
 \times 10^{-3}$   & $3.6 \times 10^{-3}$ \\
 \cline{2-4}
 & $1.5$ day &  $8.6\times 10^{-3}$   &  $5.5 \times 10^{-3}$\\
   \hline

 \multirow{2}{1.5cm}{3} & $5$ min& $4.9
 \times 10^{-3}$ & $1.6\times 10^{-3}$  \\
 \cline{2-4}
 & $26~$min& $9.6\times 10^{-3}$   &$7.4\times 10^{-3}$   \\
 \cline{2-4}
  \hline
   \multirow{2}{1.5cm}{4} & $12$ h  & $1.4\times 10^{-3}$   & $1.2\times 10^{-3}$ \\
 \cline{2-4}
   \cline{2-4}
 & $2$ days  & $3.5\times 10^{-3}$   & $1.7\times 10^{-3}$\\
  \hline
\end{tabular}
\end{center}
\end{table}

We compute the evolution of the total mass of water $I$ for each numerical solution: 
\begin{equation}
I(t)=\int_{0}^{L}\theta(z,t)dz.
\end{equation} 

\begin{figure}[ht]
\centering
\begin{tabular}{cc}
Clay & Clay loam \\
\includegraphics[width=5.6cm,height=4.2cm,angle=0]{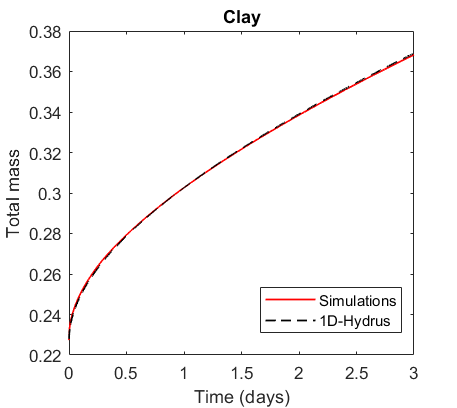} & 
\includegraphics[width=5.6cm,height=4.2cm,angle=0]{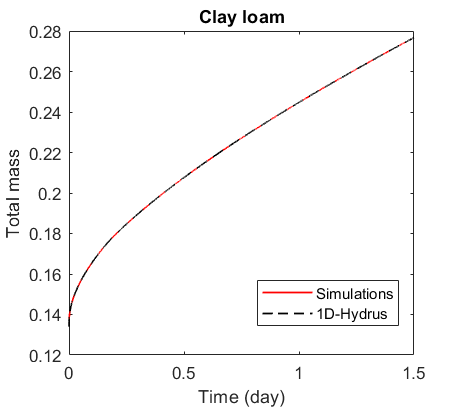} \\
Sand & Silty clay \\
\includegraphics[width=5.6cm,height=4.2cm,angle=0]{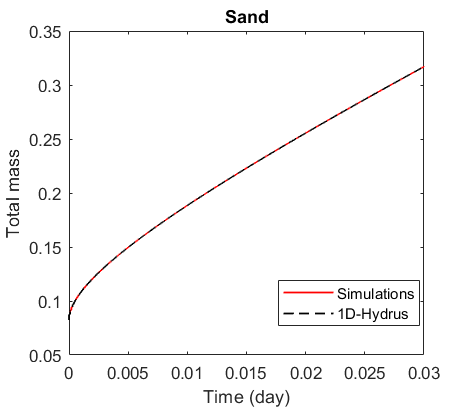} & 
\includegraphics[width=5.6cm,height=4.2cm,angle=0]{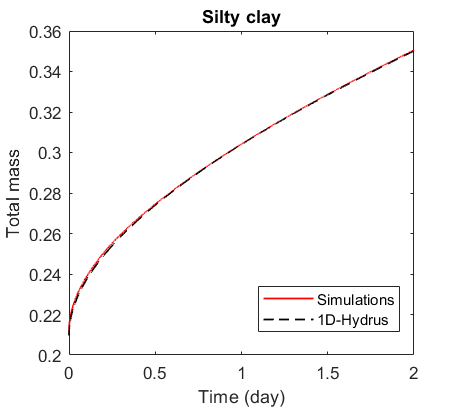} \\
\end{tabular}
\caption{Comparison of the total mass between the numerical and reference solutions.} \label{f21}
\end{figure}

\begin{figure}[ht]
\centering
\begin{tabular}{cc}
\includegraphics[width=6cm,height=4.3cm,angle=0]{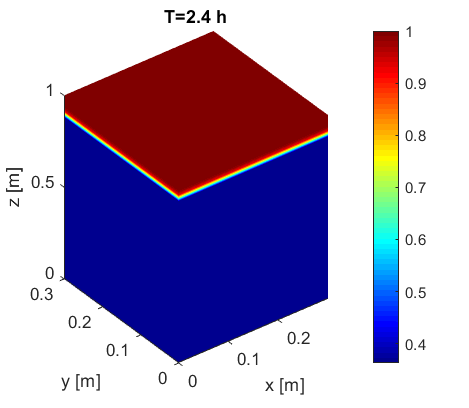} & \includegraphics[width=6cm,height=4.3cm,angle=0]{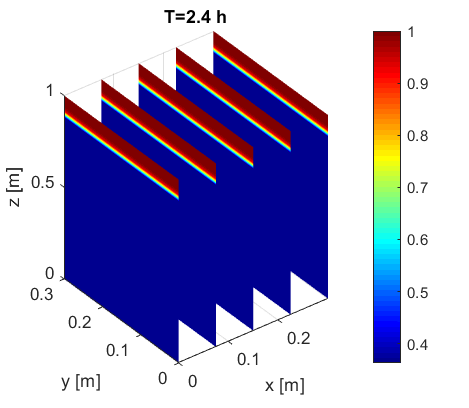}\\
\includegraphics[width=6cm,height=4.3cm,angle=0]{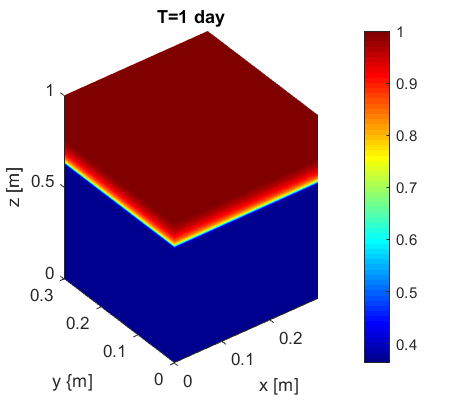} & \includegraphics[width=6cm,height=4.3cm,angle=0]{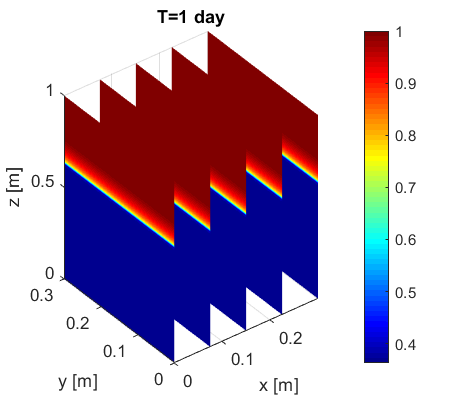}\\
\includegraphics[width=6cm,height=4.3cm,angle=0]{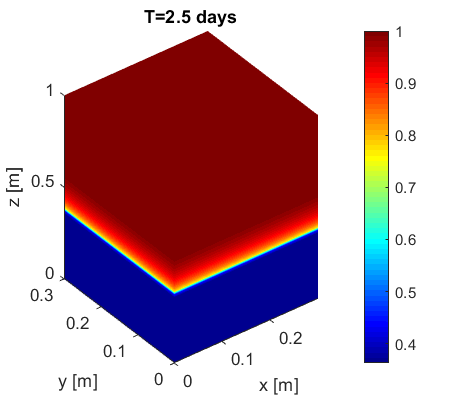} & \includegraphics[width=6cm,height=4.3cm,angle=0]{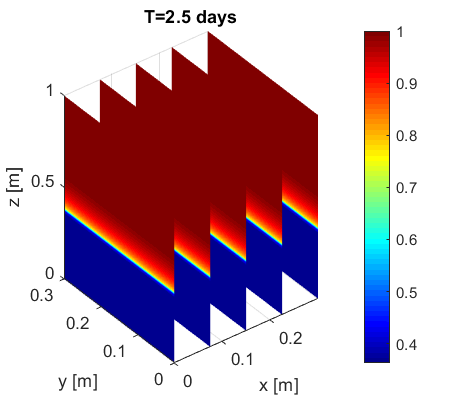}
\end{tabular}
\caption{The 3D evolution of saturation of the silty clay soil. }\label{X3D2}
\end{figure}

\begin{figure}[ht]
\centering
\begin{tabular}{cc}
\includegraphics[width=6cm,height=4.3cm,angle=0]{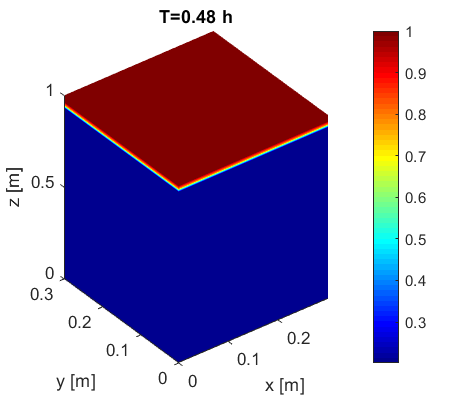} & \includegraphics[width=6cm,height=4.3cm,angle=0]{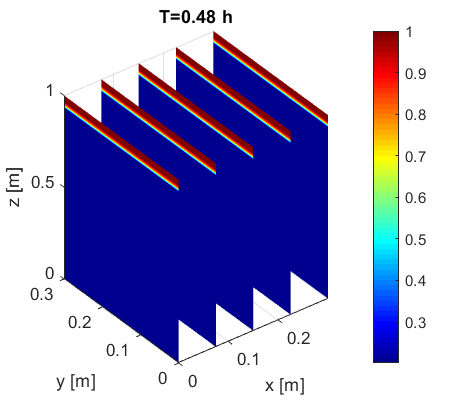}\\
\includegraphics[width=6cm,height=4.3cm,angle=0]{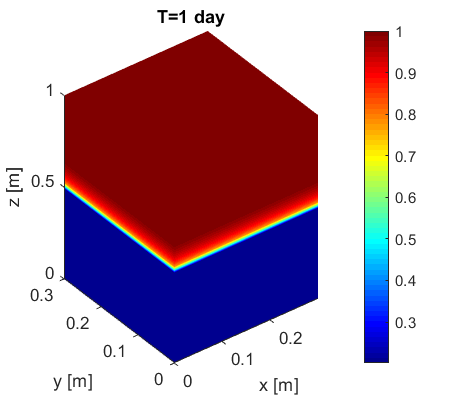} & \includegraphics[width=6cm,height=4.3cm,angle=0]{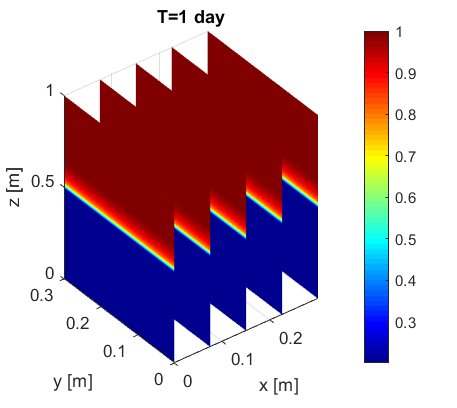}\\
\includegraphics[width=6cm,height=4.3cm,angle=0]{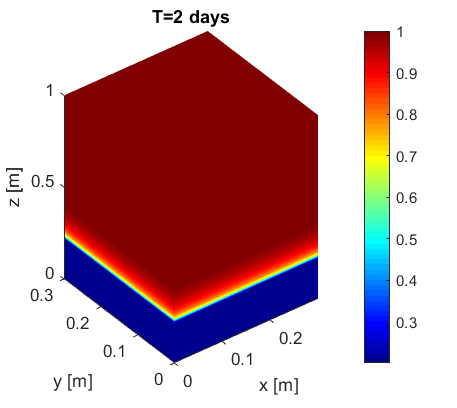} & \includegraphics[width=6cm,height=4.3cm,angle=0]{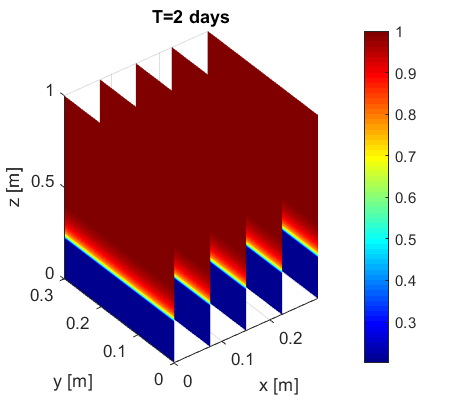}
\end{tabular}
\caption{The 3D evolution of saturation of the clay loam soil. }\label{X3D3}
\end{figure}

\begin{figure}[ht]
\centering
 \begin{tabular}{ccc}
 \\
\includegraphics[width=6cm,height=5cm,angle=0]{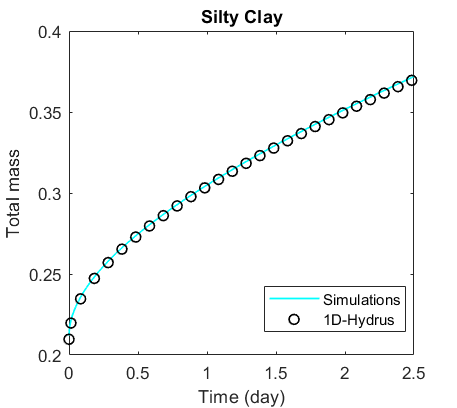} & \includegraphics[width=6cm,height=5cm,angle=0]{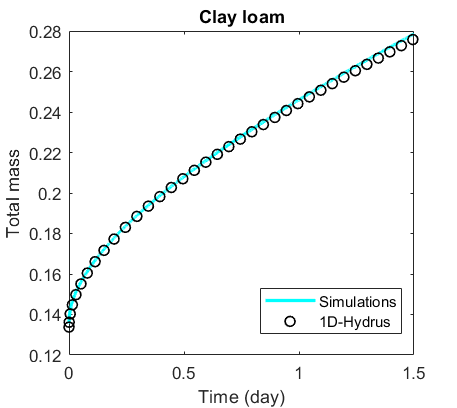}& 
\end{tabular}
\caption{Time-evolution of the total mass of water of the selected soils.}\label{Test3T}
\end{figure}

Figure \ref{f21} shows the evolution of the total mass of water for the numerical solutions obtained using the proposed techniques and $1D$-solution obtained using Hydrus.  
The results show the effectiveness of the LRBF method in terms of conservation of mass.
In the following, we perform numerical simulations for $3D$ infiltration problem. We consider a block of soil having the dimensions $l_{1}=l_{2}=0.3~m$ and $L=1~m$. We consider the same physical parameters as the previous test. The  silty clay and clay loam soils are chosen for this numerical test. We set $c=0.6$, $n_{s} =7$, $N_{x}=N_{y}=90$, $N_{z} = 300$ and $\Delta t=0.0001$.  Figures \ref{X3D2} and \ref{X3D3} show the $3D$ evolution of saturation (left) for the selected soils. The results on the right side are the $x$-slices of saturation ($x = 0$, $x = l_{1}/4$, $x = l_{1}/2$, $x = 3l_{1}/4$, $x = l_{1}$). The cross sectional average in the vertical direction of the total mass of water of $3D$ numerical solutions and the $1D$-Hydrus reference solutions ($l_x=l_y=1$) are shown in Figure \ref{Test3T}. The results confirm the accuracy of the proposed numerical method for infiltration through three-dimensional porous medium.

\subsection{Infiltration in heterogeneous soil}\label{ss42}
In the following sections, we perform numerical tests to study the robustness of the developed numerical method in modeling one-, two- and three-dimensional heterogeneous medium.
\subsubsection{Infiltration in 1D-layered soils}
\begin{figure}[ht]
\centering
\begin{tabular}{cc}
Water content ($\theta$) & Pressure head ($h$) \\
\includegraphics[width=8cm,height=6.5cm,angle=0]{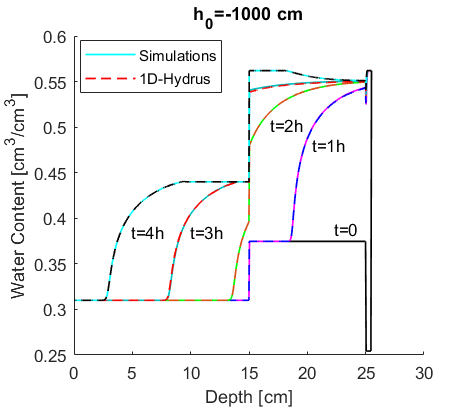} &\includegraphics[width=8cm,height=6.5cm,angle=0]{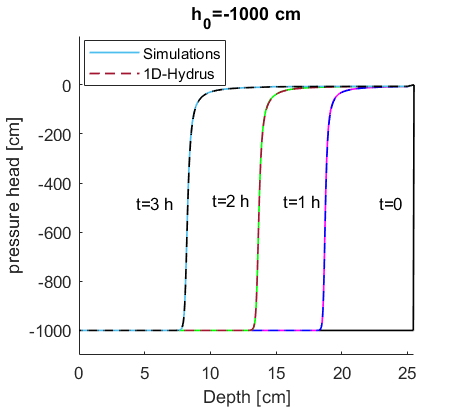}  \\
\includegraphics[width=8cm,height=6.5cm,angle=0]{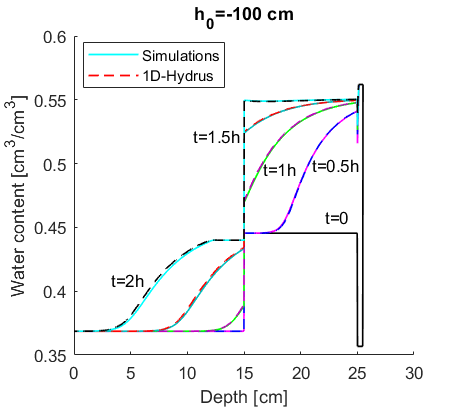} &\includegraphics[width=8cm,height=6.5cm,angle=0]{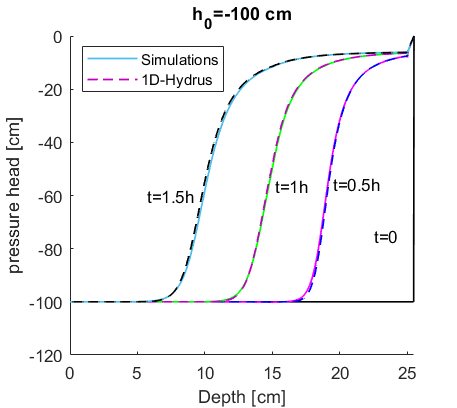}

\end{tabular}
\caption{Vertical profile of soil water content (left) and pressure head (right).}\label{1Dh}

\end{figure}

In this numerical test, we perform simulations of infiltration using a column of soil ($L=25.5~cm$) with three layers. The layered soil consists of a thin surface crust ($0.5~cm$), a tilled layer ($10~cm$) and a subsoil layer ($15~cm$). The hydraulic properties \cite{manzini2004mass} for the layers soil are shown in Table \ref{t3n}.
\begin{table}[ht]
\begin{center}
\caption{Parameters of the  layered soil. }\label{t3n}
 \begin{tabular}{|c|c|c|c|c|c|c|}
 \hline
Layer & Elevation  & $\theta_s$ & $K_s$  & $h_{\text{d}}$  & $\lambda$ & $\beta$ \\
 --& $(cm)$ & --  & $(cm/h)$ & $(cm)$ & --  & -- \\ 
 \hline
Surface crust & $25 \leq z \leq 25.5$ & $0.562$ & $0.0616$ & $-4.55$ & $0.1470$ & $ 16.6054$\\
 \hline
Tilled layer & $15 \leq z \leq 25$ & $0.562$ & $1.396$ & $-4.55$ & $0.0751$ & $29.6312$  \\
 \hline
Sub-soil &  $0 \leq z \leq 15$& $0.440$ & $0.312$ & $-9.50$ & $0.0751$ & $29.6312$  \\
 \hline
  \end{tabular}
\end{center}
 \end{table}
Numerical simulations are performed for two cases using $h_{0}=-100~cm$ and $-1000~cm$. We set $c=0.6$, $n_{s} =3$, $N_{z} = 1001$ and $\Delta t=0.005$. Figure \ref{1Dh} displays the evolution in time of the water content (left) and pressure head (right). The numerical solutions are in good agreement with the $1D$-Hydrus simulations.
Table \ref{RMSE2} presents the $RMSE$ and $L^{1}_{er}$ errors between the numerical solutions and the results obtained using $1D$-Hydrus for the soil water content. The results confirm the effectiveness of the developed method in terms of accuracy in modeling infiltration in layered soils.
\begin{table}[ht]
\begin{center}
\caption{The computed errors between the numerical and reference solutions.}\label{RMSE2}
 \begin{tabular}{|c|c|c|c|}
 \hline
 $h_{0}$ & $T$ & $RMSE$ & $L^{1}_{er}$   \\
 \hline
 \multirow{2}{1.5cm}{$-100~m $} & $0.5~h$& $4.72\times 10^{-4}$ & $4.06\times 10^{-4}$  \\
 \cline{2-4}
 & $1~h$& $6.56\times 10^{-4}$  &$7.34\times 10^{-4}$  \\
 \cline{2-4}
 & $1.5~h$& $9.96\times 10^{-4}$    &$1.3\times 10^{-3}$   \\
 \hline
 \multirow{2}{1.5cm}{$-1000~m $} & $1~h$  & $3.3\times 10^{-3}$   & $2.4\times 10^{-3}$  \\
 \cline{2-4}
 & $2~h$  & $1.2\times 10^{-3}$   & $1.1\times 10^{-3}$  \\

   \cline{2-4}
 & $3~h$  & $1.5\times 10^{-3}$   & $1.9\times 10^{-3}$ \\
  \hline
\end{tabular}
\end{center}
\end{table}
\subsubsection{ Infiltration in 2D-layered soils}
\begin{figure}[ht]
\centering
\begin{tabular}{ccc}
\includegraphics[width=5cm,height=4cm,angle=0]{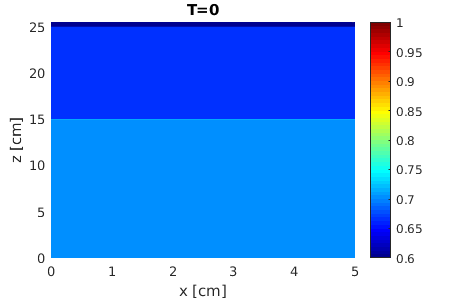} & \includegraphics[width=5cm,height=4cm,angle=0]{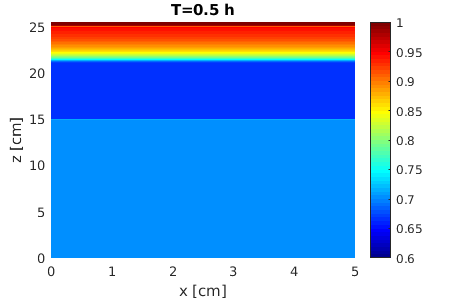} &
\includegraphics[width=5cm,height=4cm,angle=0]{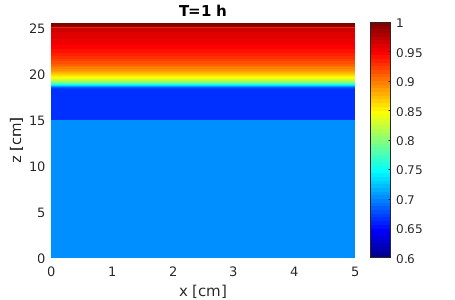} \\ \includegraphics[width=5cm,height=4cm,angle=0]{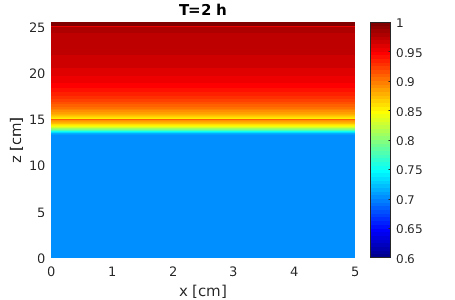}&
\includegraphics[width=5cm,height=4cm,angle=0]{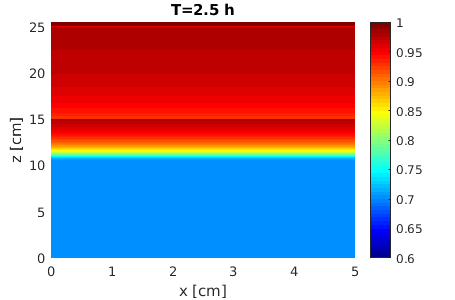} & \includegraphics[width=5cm,height=4cm,angle=0]{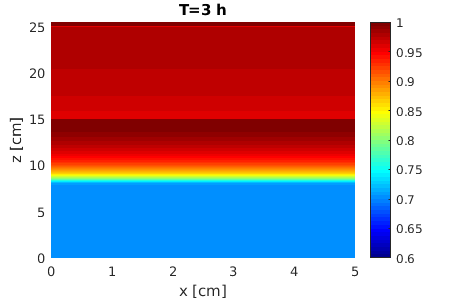}
\end{tabular}
\caption{The time evolution of saturation of the considered soil.}\label{X2Dh}
\end{figure}

This numerical test is performed for unsaturated flow through $2D$ layered porous medium. We consider the physical parameters of soils given in Table \ref{t3n} and $l_{1}=5~cm$. Figure \ref{X2Dh} displays the time-evolution of saturation for $h_{0}=-1000~cm$ obtained using the proposed method. The results are obtained using $N_{x}=100$, $N_{z}=1001$, $\Delta t=0.005$, $n_{s}=5$ and $c=0.6$.
 \begin{figure}[ht]
\centering
 \begin{tabular}{ccc}
 \\
\includegraphics[width=6cm,height=4.8cm,angle=0]{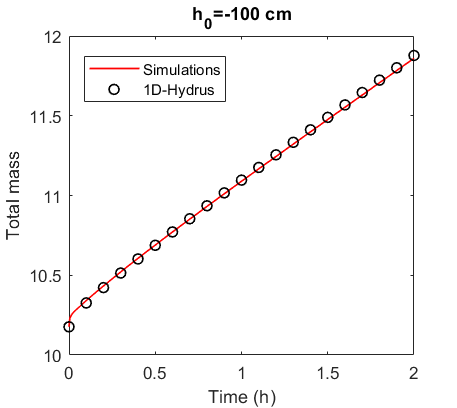} & \includegraphics[width=6cm,height=4.8cm,angle=0]{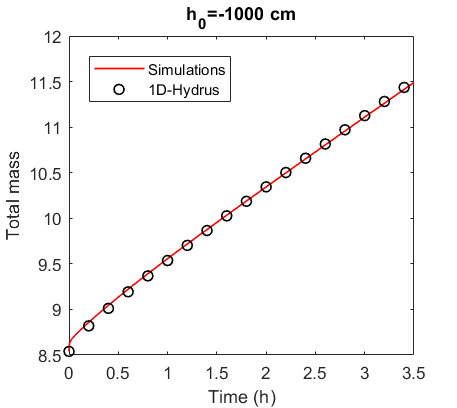}& 
\end{tabular}
\caption{Time-evolution of the total mass of water for $h_{0}=-100~cm$ and $-1000~cm$.}\label{TM21}
\end{figure}
We present in Figure \ref{TM21} the cross sectional average in the vertical direction of the total mass of water of $2D$ numerical solutions and the $1D$-Hydrus reference solutions ($l_x=1$) where we consider two cases using $h_{0}=-100~cm$ and $-1000~cm$. The results show a good correspondence between the numerical and reference solutions which confirms the accuracy of the proposed method.
\subsubsection{Infiltration in 3D-layered soils}
\begin{figure}[ht]
\centering
\begin{tabular}{cc}
\includegraphics[width=6cm,height=4.3cm,angle=0]{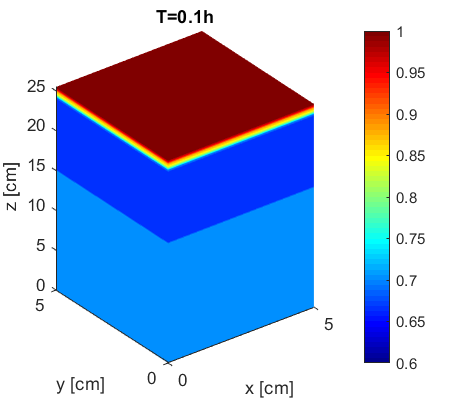} & \includegraphics[width=6cm,height=4.3cm,angle=0]{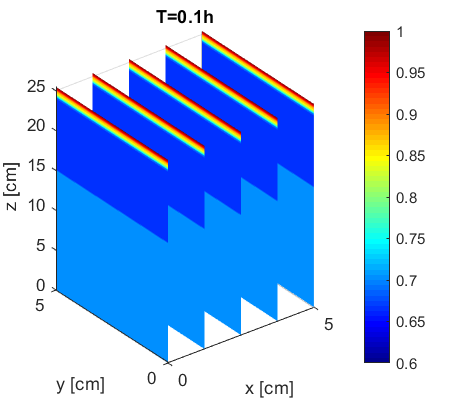}\\
\includegraphics[width=6cm,height=4.3cm,angle=0]{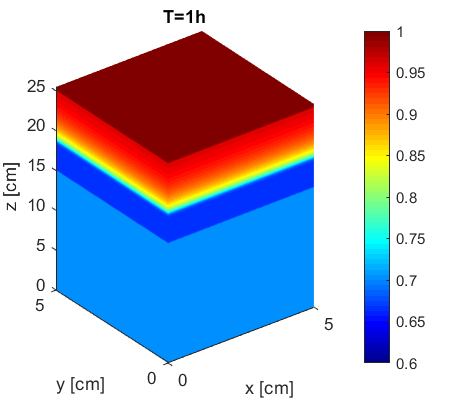} & \includegraphics[width=6cm,height=4.3cm,angle=0]{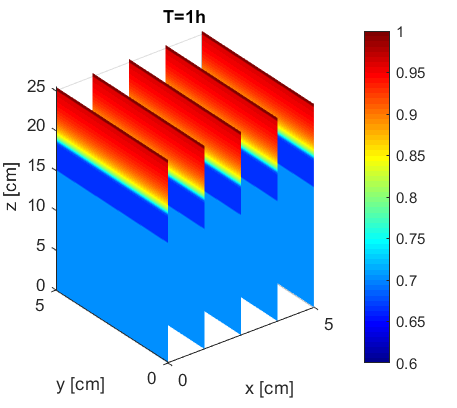}\\
\includegraphics[width=6cm,height=4.3cm,angle=0]{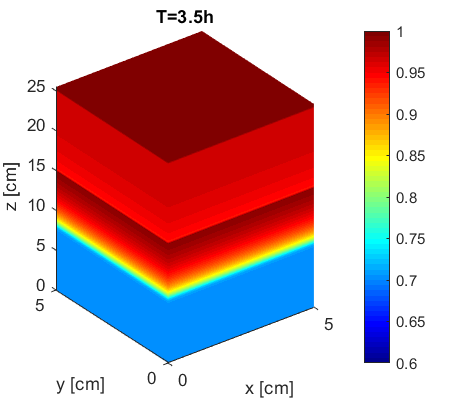} & \includegraphics[width=6cm,height=4.3cm,angle=0]{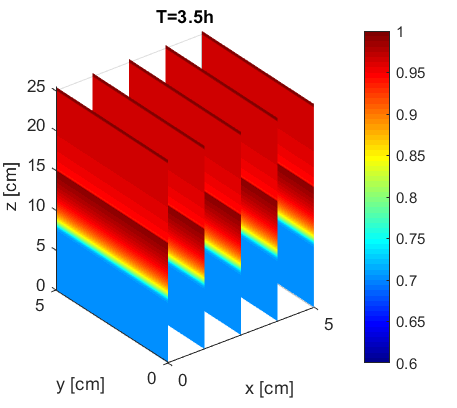}
\end{tabular}
\caption{The 3D evolution of saturation of the layered soils. }\label{A3D2}
\end{figure}

\begin{figure}[ht]
\centering

\includegraphics[width=6cm,height=4.8cm,angle=0]{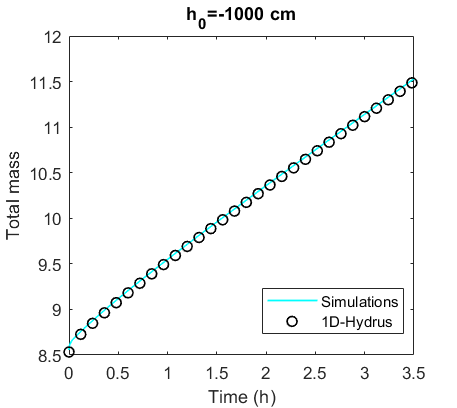} 
\caption{Time-evolution of the total mass of water for $h_{0}=-1000~cm$.}\label{T3}
\end{figure}

Here, we investigate the capability of the developed numerical model in predicting infiltration through three-dimensional layered porous medium. We consider the same hydraulic properties of soils as the previous test. We perform numerical simulations using $c=0.6$, $n_{s} =7$, $N_{x}=N_{y}=100$, $N_{z} = 501$ and $\Delta t=0.001$. In Figure \ref{A3D2}, we display the $3D$ evolution of saturation (left) for the considered soils for $h_{0}=-1000~cm$. The results on the right side present the $x$-slices of saturation ($x = 0~$, $x = l_{1}/4$, $x = l_{1}/2$, $x = 3l_{1}/4$, $x = l_{1}$).
The cross sectional average in the vertical direction of the total mass of water of $3D$ numerical solutions and the $1D$-Hydrus reference solutions ($l_x=l_y=1$) are shown in Figure \ref{T3}. 
The comparison between the results of the total mass shows the accuracy of the developed numerical model for infiltration in three-dimensional layered soils.
\subsubsection{Infiltration in layered soil of L-shape form}
\begin{figure}[ht]
\centering

\includegraphics[width=6cm,height=6cm,angle=0]{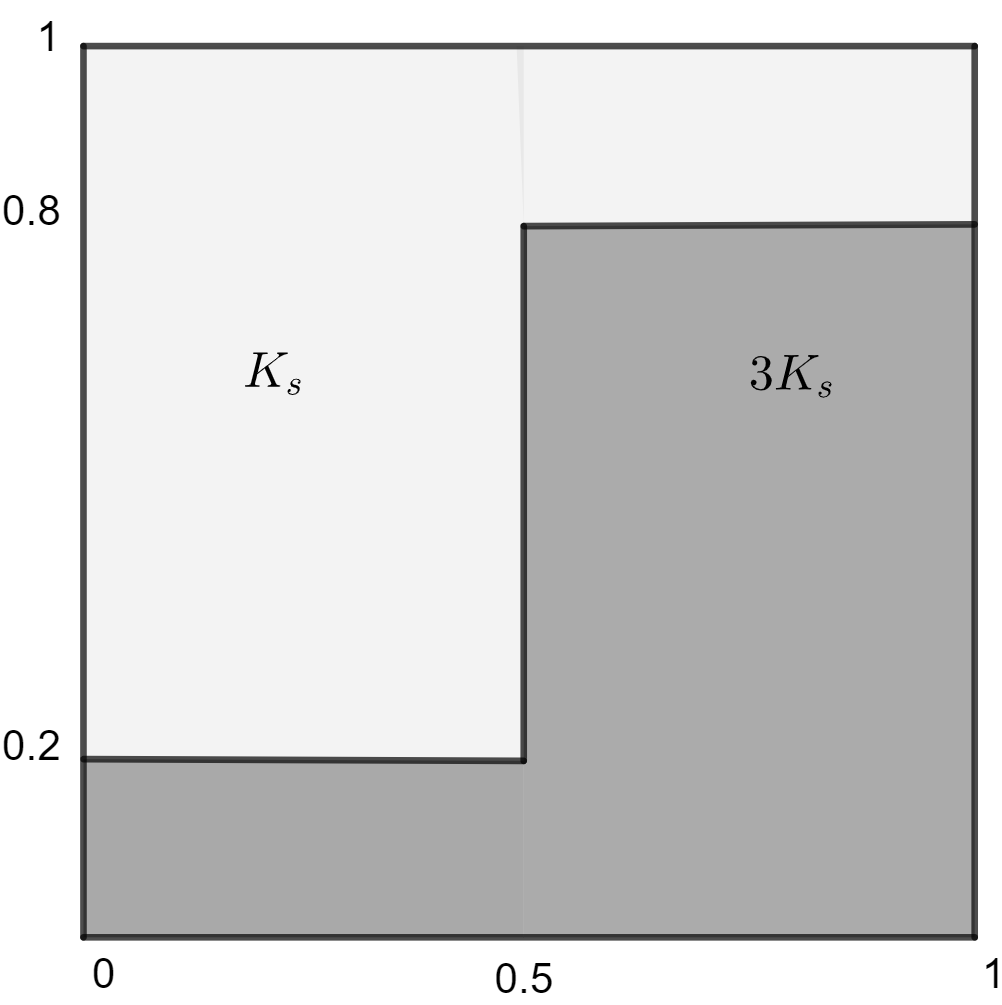} 
\caption{Domain description.}\label{Lshape}
\end{figure}

Here, we used a complex geometry of the interface between soils compared to the previous test. We perform numerical simulations using infiltration problem through layered soil of $L$-shape \cite{baron2017adaptive}. As shown in Figure \ref{Lshape}, the computational domain is partitioned into two subdomains  which differ in their saturated hydraulic conductivity $K_s$. 
In \cite{baron2017adaptive}, the van Genuchten model \cite{van1980closed} is used for capillary pressure with the parameter values $K_s=  0.3319~m/h$, $\theta_s=0.368$, $\theta_r=0.102$, $\alpha=3.35 ~m^{-1}$, $n=2$ and $m=0.5$, where $\alpha=1/h_{\text{cap}}$, $n$ and $m$ are empirical parameters satisfy $m=1-1/n$.
 In our case, we used the Brooks-Corey model \cite{brooks1964hydraulic} where we approximate  $\lambda$ and $h_d$ based on the equivalence between  van Genuchten and Brooks-Corey parameters proposed in \cite{lenhard1989corey}. The parameters $h_d$ and $\lambda$ are given by \cite{lenhard1989corey}:
\begin{equation}\label{BCVGhd}
    h_d=\left(\dfrac{1}{\alpha}\right)S_x^{1/\lambda}(S_x^{-1/m}-1 )^{1-m},
\end{equation}
\begin{equation}\label{BCVGlmb}
    \lambda=\dfrac{m}{1-m}(1-0.5^{1/m}),
\end{equation}
where $S_x=0.72-0.35\exp(-n^4)$. We used a homogeneous Neumann condition on the vertical sides of the domain ($x=0, 1~m$) and a homogeneous Dirichlet on the top and bottom sides ($z=0, 1~m$). The initial condition is $h(x,z,0)=-z$. We set $\varepsilon=0.1$, $n_s = 5$, $N_x=N_z=1000$ and $\Delta t = 0.001$. Figure \ref{LsRes} displays the evolution of saturation at times $T=12,~24,\text{and}~48~h$. We observe good overall agreement between the results of our simulations and those presented in 
 \cite{baron2017adaptive,keita2021implicit}.
 \begin{figure}[ht]
\centering
\begin{tabular}{ccc}
$t=12~h$ & $t=24~h$ & $t=48~h$ \\
    \includegraphics[width=6cm,height=5cm,angle=0]{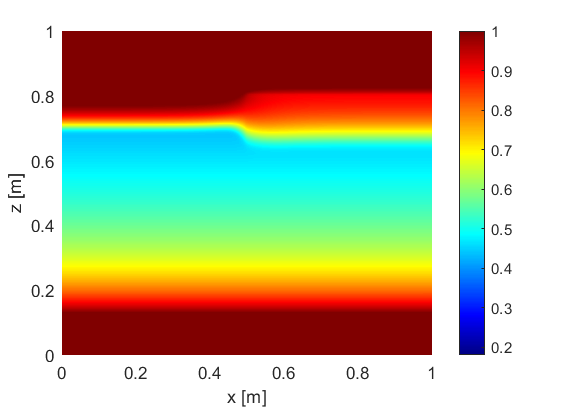} &
\includegraphics[width=6cm,height=5cm,angle=0]{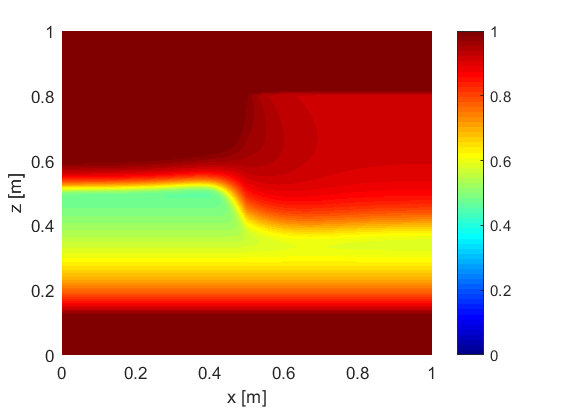} & \includegraphics[width=6cm,height=5cm,angle=0]{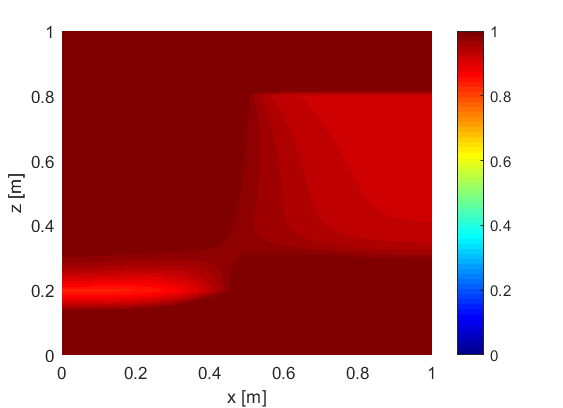} 

\end{tabular}
\caption{Water saturation at $T=12,24,48~h$. }\label{LsRes}
\end{figure}
\subsubsection{Infiltration in curvilinearly layered soil}
\begin{figure}[ht]
\centering

\includegraphics[width=6cm,height=6cm,angle=0]{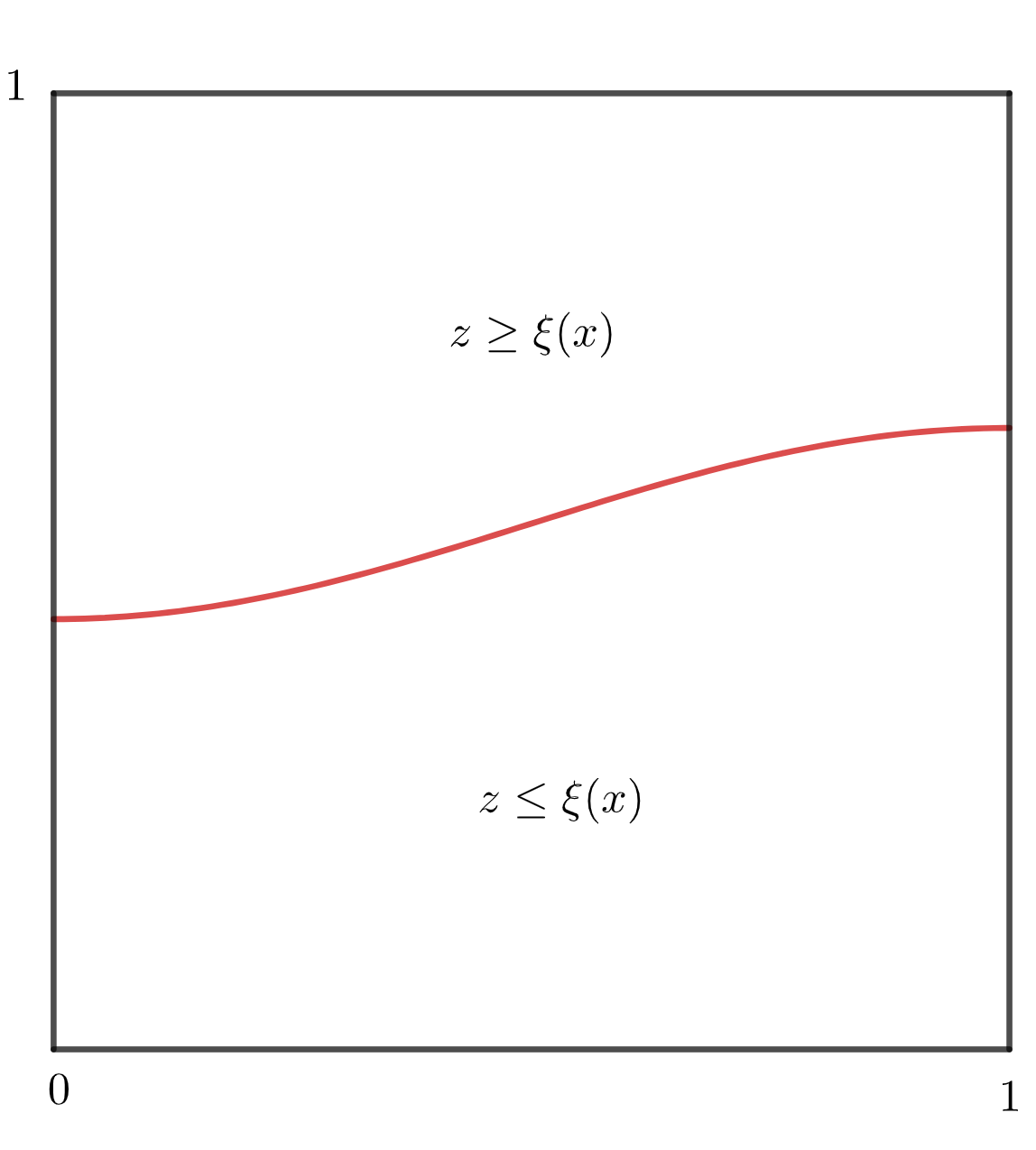}
\caption{Domain description.}\label{CurvDom}
\end{figure}

In this last numerical test, we perform simulations of infiltration through curvilinearly layered soil. The computational domain is split into two subdomains separated by a curved interface (see Figure \ref{CurvDom}).
The interface equation is given by:
\begin{equation}
    \xi(x)=l_2\left(0.1\left(1-\cos(\pi x/l_1)\right)+0.45\right).
\end{equation}
The hydraulic parameters of the soils \cite{keita2021implicit} are shown in Table \ref{t4}.
\begin{table}[ht]
\begin{center}
\caption{Parameters of the soils. }\label{t4}
 \begin{tabular}{|c|c|c|c|c|c|c|}
 \hline
 Elevation  & $\theta_r$ & $\theta_s$ & $K_s$  & $h_{\text{d}}$  & $\lambda$ & $\beta$  \\
 $(m)$& -- &  --   & $(m/h)$ & $(m)$ & --  & -- \\ 
 \hline
 $ z \geq \xi(x)$ &$0120$ & $0.5$ & $0.0025$ & $-0.45$ & $0.34$ & $ 3.02$\\
 \hline
 $ z \leq \xi(x)$ & $0.034$   & $0.46$ & $0.02$ & $-0.23$ & $1.29$ & $5.87$ \\
 \hline
  \end{tabular}
\end{center}
 \end{table}
 As in the previous test, the values of the parameters $h_d$ and $\lambda$ are approximated using Equations \eqref{BCVGhd} and \eqref{BCVGlmb}. We use $l_1=l_2=1~m$, $\varepsilon=0.8$, $n_s = 5$, $N_x=N_z=1000$ and $\Delta t = 0.001$. We consider the initial condition $h(x,z,0)=-z$. Homogeneous Dirichlet boundary conditions are imposed on
the top and bottom sides of the domain, while  homogeneous Neumann conditions are enforced on the two vertical sides. We present in Figure \ref{CurvDom} the time evolution of the water content.
  \begin{figure}[ht]
\centering
\begin{tabular}{ccc}
$t=0$ & $t=0.25$~day & $t=0.75$~day  \\
    \includegraphics[width=5.5cm,height=4cm,angle=0]{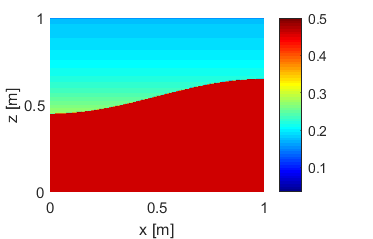} &
\includegraphics[width=5.5cm,height=4cm,angle=0]{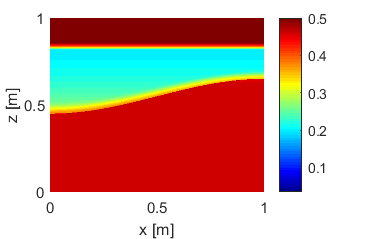} &  \includegraphics[width=5.5cm,height=4cm,angle=0]{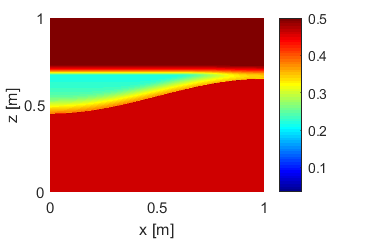} \\ 
$t=1$~day & $t=1.25$~day   & $t=1.75$~day  \\
    \includegraphics[width=5.5cm,height=4cm,angle=0]{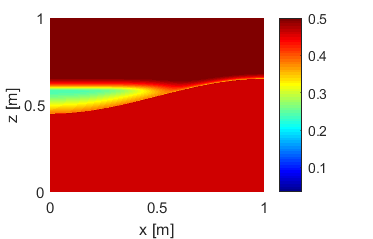} &
\includegraphics[width=5.5cm,height=4cm,angle=0]{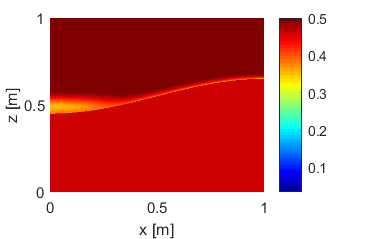} &  \includegraphics[width=5.5cm,height=4cm,angle=0]{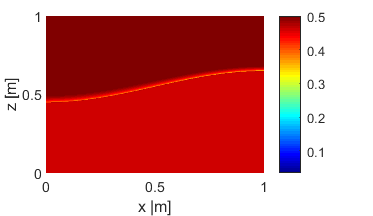} 

\end{tabular}
\caption{Time evolution of water content. }\label{LsRes}
\end{figure}
The results of our simulations are in good overall agreement at different times compared to the results of simulation presented in \cite{manzini2004mass, keita2021implicit}, which confirm the robustness of the developed numerical model in predicting infiltration in heterogeneous soils.
\section{Conclusion}\label{sec:5}
In this study, a new approach is developed for modeling unsaturated flow through porous media. The proposed techniques are based on the Kirchhoff transformation, the Brooks and Corey model for the capillary pressure function and a power-law relation for the relative permeability function. The proposed approach allows us to avoid technical issues associated with the use of the Kirchhoff transformation in heterogeneous soils and to reduce the nonlinearity of the model equation. The resulting system is solved based on the LRBF methods which are very effective for solving high-dimensional problems since they don't require mesh generation and have a computational advantage of using reduced memory. The LRBF meshless methods allow us to avoid ill-conditioning problems where a sparse matrix is obtained for the global system. The performance and robustness of the developed numerical model are demonstrated based on comparisons between numerical and reference solutions. Numerical experiments are performed to simulate the infiltration in one, two, and
three-dimensional soils. The numerical results show the accuracy of the proposed techniques for modeling infiltration through porous media.

\section*{Acknowledgment}
AB gratefully acknowledges funding from UM6P-OCP.
\bibliographystyle{elsarticle-num}
\bibliography{mabiblio.bib}

\begin{thebibliography}{10}
\expandafter\ifx\csname url\endcsname\relax
  \def\url#1{\texttt{#1}}\fi
\expandafter\ifx\csname urlprefix\endcsname\relax\def\urlprefix{URL }\fi
\expandafter\ifx\csname href\endcsname\relax
  \def\href#1#2{#2} \def\path#1{#1}\fi

\bibitem{richards1931capillary}
L.~A. Richards, Capillary conduction of liquids through porous mediums, Physics
  1~(5) (1931) 318--333.
\newblock \href {https://doi.org/https://doi.org/10.1063/1.1745010}
  {\path{doi:https://doi.org/10.1063/1.1745010}}.

\bibitem{gardner1958some}
W.~Gardner, Some steady-state solutions of the unsaturated moisture flow
  equation with application to evaporation from a water table, Soil science
  85~(4) (1958) 228--232.

\bibitem{brooks1964hydraulic}
R.~Brooks, A.~Corey, Hydraulic properties of porous media: Hydrology papers:
  Colorado state university, Fort Collins, Colorado (1964).

\bibitem{van1980closed}
M.~T. Van~Genuchten, A closed-form equation for predicting the hydraulic
  conductivity of unsaturated soils, Soil science society of America journal
  44~(5) (1980) 892--898.
\newblock \href
  {https://doi.org/https://doi.org/10.2136/sssaj1980.03615995004400050002x}
  {\path{doi:https://doi.org/10.2136/sssaj1980.03615995004400050002x}}.

\bibitem{SrivastavaYeh}
R.~Srivastava, T.~J. Yeh, Analytical solutions for one‐dimensional, transient
  infiltration toward the water table in homogeneous and layered soils, Water
  Resources Research 27~(5) (1991) 753--762.
\newblock \href {https://doi.org/https://doi.org/10.1029/90WR02772}
  {\path{doi:https://doi.org/10.1029/90WR02772}}.

\bibitem{tracy19951}
F.~T. Tracy, 1-{D}, 2-{D}, and 3-{D} analytical solutions of unsaturated flow
  in groundwater, Journal of hydrology 170~(1-4) (1995) 199--214.
\newblock \href {https://doi.org/https://doi.org/10.1016/0022-1694(94)02674-Z}
  {\path{doi:https://doi.org/10.1016/0022-1694(94)02674-Z}}.

\bibitem{huang2012analytical}
R.~Huang, L.~Wu, Analytical solutions to 1-{D} horizontal and vertical water
  infiltration in saturated/unsaturated soils considering time-varying
  rainfall, Computers and Geotechnics 39 (2012) 66--72.
\newblock \href {https://doi.org/https://doi.org/10.1016/j.compgeo.2011.08.008}
  {\path{doi:https://doi.org/10.1016/j.compgeo.2011.08.008}}.

\bibitem{hayek2016exact}
M.~Hayek, An exact explicit solution for one-dimensional, transient, nonlinear
  {R}ichards’ equation for modeling infiltration with special hydraulic
  functions, Journal of Hydrology 535 (2016) 662--670.
\newblock \href {https://doi.org/https://doi.org/10.1016/j.jhydrol.2016.02.021}
  {\path{doi:https://doi.org/10.1016/j.jhydrol.2016.02.021}}.

\bibitem{rucker2005parameter}
D.~F. Rucker, A.~W. Warrick, T.~P. Ferr{\'e}, Parameter equivalence for the
  {G}ardner and van {G}enuchten soil hydraulic conductivity functions for
  steady vertical flow with inclusions, Advances in water resources 28~(7)
  (2005) 689--699.
\newblock \href
  {https://doi.org/https://doi.org/10.1016/j.advwatres.2005.01.004}
  {\path{doi:https://doi.org/10.1016/j.advwatres.2005.01.004}}.

\bibitem{haverkamp1977comparison}
R.~Haverkamp, M.~Vauclin, J.~Touma, P.~Wierenga, G.~Vachaud, A comparison of
  numerical simulation models for one-dimensional infiltration, Soil Science
  Society of America Journal 41~(2) (1977) 285--294.
\newblock \href
  {https://doi.org/https://doi.org/10.2136/sssaj1977.03615995004100020024x}
  {\path{doi:https://doi.org/10.2136/sssaj1977.03615995004100020024x}}.

\bibitem{celia1990general}
M.~A. Celia, E.~T. Bouloutas, R.~L. Zarba, A general mass-conservative
  numerical solution for the unsaturated flow equation, Water resources
  research 26~(7) (1990) 1483--1496.
\newblock \href {https://doi.org/https://doi.org/10.1029/WR026i007p01483}
  {\path{doi:https://doi.org/10.1029/WR026i007p01483}}.

\bibitem{clement1994physically}
T.~Clement, W.~R. Wise, F.~J. Molz, A physically based, two-dimensional,
  finite-difference algorithm for modeling variably saturated flow, Journal of
  Hydrology 161~(1-4) (1994) 71--90.
\newblock \href {https://doi.org/https://doi.org/10.1016/0022-1694(94)90121-X}
  {\path{doi:https://doi.org/10.1016/0022-1694(94)90121-X}}.

\bibitem{huyakorn1984techniques}
P.~Huyakorn, S.~Thomas, B.~Thompson, Techniques for making finite elements
  competitve in modeling flow in variably saturated porous media, Water
  Resources Research 20~(8) (1984) 1099--1115.
\newblock \href {https://doi.org/https://doi.org/10.1029/WR020i008p01099}
  {\path{doi:https://doi.org/10.1029/WR020i008p01099}}.

\bibitem{radu2004order}
F.~Radu, I.~S. Pop, P.~Knabner, Order of convergence estimates for an {E}uler
  implicit, mixed finite element discretization of {R}ichards' equation, SIAM
  Journal on Numerical Analysis 42~(4) (2004) 1452--1478.
\newblock \href {https://doi.org/https://doi.org/10.1137/S0036142902405229}
  {\path{doi:https://doi.org/10.1137/S0036142902405229}}.

\bibitem{bause2004computation}
M.~Bause, P.~Knabner, Computation of variably saturated subsurface flow by
  adaptive mixed hybrid finite element methods, Advances in Water Resources
  27~(6) (2004) 565--581.
\newblock \href
  {https://doi.org/https://doi.org/10.1016/j.advwatres.2004.03.005}
  {\path{doi:https://doi.org/10.1016/j.advwatres.2004.03.005}}.

\bibitem{solin2011solving}
P.~Solin, M.~Kuraz, Solving the nonstationary {R}ichards equation with adaptive
  hp-{FEM}, Advances in water resources 34~(9) (2011) 1062--1081.
\newblock \href
  {https://doi.org/https://doi.org/10.1016/j.advwatres.2011.04.020}
  {\path{doi:https://doi.org/10.1016/j.advwatres.2011.04.020}}.

\bibitem{eymard1999finite}
R.~Eymard, M.~Gutnic, D.~Hilhorst, The finite volume method for {R}ichards
  equation, Computational Geosciences 3~(3) (1999) 259--294.
\newblock \href {https://doi.org/https://doi.org/10.1023/A:1011547513583}
  {\path{doi:https://doi.org/10.1023/A:1011547513583}}.

\bibitem{manzini2004mass}
G.~Manzini, S.~Ferraris, Mass-conservative finite volume methods on 2-{D}
  unstructured grids for the {R}ichards’ equation, Advances in Water
  Resources 27~(12) (2004) 1199--1215.
\newblock \href
  {https://doi.org/https://doi.org/10.1016/j.advwatres.2004.08.008}
  {\path{doi:https://doi.org/10.1016/j.advwatres.2004.08.008}}.

\bibitem{lai2015mass}
W.~Lai, F.~L. Ogden, A mass-conservative finite volume predictor--corrector
  solution of the 1{D} {R}ichards’ equation, Journal of Hydrology 523 (2015)
  119--127.
\newblock \href {https://doi.org/https://doi.org/10.1016/j.jhydrol.2015.01.053}
  {\path{doi:https://doi.org/10.1016/j.jhydrol.2015.01.053}}.

\bibitem{svyatskiy2017second}
D.~Svyatskiy, K.~Lipnikov, Second-order accurate finite volume schemes with the
  discrete maximum principle for solving {R}ichards’ equation on unstructured
  meshes, Advances in water resources 104 (2017) 114--126.
\newblock \href
  {https://doi.org/https://doi.org/10.1016/j.advwatres.2017.03.015}
  {\path{doi:https://doi.org/10.1016/j.advwatres.2017.03.015}}.

\bibitem{ngo2020control}
D.~Ngo-Cong, N.~Mai-Duy, D.~L. Antille, M.~T. van Genuchten, A control volume
  scheme using compact integrated radial basis function stencils for solving
  the {R}ichards equation, Journal of Hydrology 580 (2020) 124240.
\newblock \href {https://doi.org/https://doi.org/10.1016/j.jhydrol.2019.124240}
  {\path{doi:https://doi.org/10.1016/j.jhydrol.2019.124240}}.

\bibitem{kirkland1992algorithms}
M.~R. Kirkland, R.~Hills, P.~Wierenga, Algorithms for solving {R}ichards'
  equation for variably saturated soils, Water Resources Research 28~(8) (1992)
  2049--2058.
\newblock \href {https://doi.org/https://doi.org/10.1029/92WR00802}
  {\path{doi:https://doi.org/10.1029/92WR00802}}.

\bibitem{huang1996new}
K.~Huang, B.~Mohanty, M.~T. Van~Genuchten, A new convergence criterion for the
  modified {P}icard iteration method to solve the variably saturated flow
  equation, Journal of Hydrology 178~(1-4) (1996) 69--91.
\newblock \href {https://doi.org/https://doi.org/10.1016/0022-1694(95)02799-8}
  {\path{doi:https://doi.org/10.1016/0022-1694(95)02799-8}}.

\bibitem{lehmann1998comparison}
F.~Lehmann, P.~Ackerer, Comparison of iterative methods for improved solutions
  of the fluid flow equation in partially saturated porous media, Transport in
  porous media 31~(3) (1998) 275--292.
\newblock \href {https://doi.org/https://doi.org/10.1023/A:1006555107450}
  {\path{doi:https://doi.org/10.1023/A:1006555107450}}.

\bibitem{an2012comparison}
H.~An, Y.~Ichikawa, Y.~Tachikawa, M.~Shiiba, Comparison between iteration
  schemes for three-dimensional coordinate-transformed saturated--unsaturated
  flow model, Journal of Hydrology 470 (2012) 212--226.
\newblock \href {https://doi.org/https://doi.org/10.1016/j.jhydrol.2012.08.056}
  {\path{doi:https://doi.org/10.1016/j.jhydrol.2012.08.056}}.

\bibitem{zha2017modified}
Y.~Zha, J.~Yang, L.~Yin, Y.~Zhang, W.~Zeng, L.~Shi, A modified {P}icard
  iteration scheme for overcoming numerical difficulties of simulating
  infiltration into dry soil, Journal of hydrology 551 (2017) 56--69.
\newblock \href {https://doi.org/https://doi.org/10.1016/j.jhydrol.2017.05.053}
  {\path{doi:https://doi.org/10.1016/j.jhydrol.2017.05.053}}.

\bibitem{ji2008generalized}
S.-H. Ji, Y.-J. Park, E.~A. Sudicky, J.~F. Sykes, A generalized transformation
  approach for simulating steady-state variably-saturated subsurface flow,
  Advances in Water Resources 31~(2) (2008) 313--323.
\newblock \href
  {https://doi.org/https://doi.org/10.1016/j.advwatres.2007.08.010}
  {\path{doi:https://doi.org/10.1016/j.advwatres.2007.08.010}}.

\bibitem{list2016study}
F.~List, F.~A. Radu, A study on iterative methods for solving {R}ichards’
  equation, Computational Geosciences 20~(2) (2016) 341--353.
\newblock \href {https://doi.org/https://doi.org/10.1007/s10596-016-9566-3}
  {\path{doi:https://doi.org/10.1007/s10596-016-9566-3}}.

\bibitem{pop2002error}
I.~S. Pop, Error estimates for a time discretization method for the {R}ichards'
  equation, Computational geosciences 6~(2) (2002) 141--160.
\newblock \href {https://doi.org/https://doi.org/10.1023/A:1019936917350}
  {\path{doi:https://doi.org/10.1023/A:1019936917350}}.

\bibitem{berninger2011fast}
H.~Berninger, R.~Kornhuber, O.~Sander, Fast and robust numerical solution of
  the {R}ichards equation in homogeneous soil, SIAM Journal on Numerical
  Analysis 49~(6) (2011) 2576--2597.
\newblock \href {https://doi.org/https://doi.org/10.1137/100782887}
  {\path{doi:https://doi.org/10.1137/100782887}}.

\bibitem{suk2019numerical}
H.~Suk, E.~Park, Numerical solution of the {K}irchhoff-transformed {R}ichards
  equation for simulating variably saturated flow in heterogeneous layered
  porous media, Journal of Hydrology 579 (2019) 124213.
\newblock \href {https://doi.org/https://doi.org/10.1016/j.jhydrol.2019.124213}
  {\path{doi:https://doi.org/10.1016/j.jhydrol.2019.124213}}.

\bibitem{ross1990efficient}
P.~J. Ross, Efficient numerical methods for infiltration using {R}ichards'
  equation, Water Resources Research 26~(2) (1990) 279--290.
\newblock \href {https://doi.org/https://doi.org/10.1029/WR026i002p00279}
  {\path{doi:https://doi.org/10.1029/WR026i002p00279}}.

\bibitem{stevens2011scalable}
D.~Stevens, H.~Power, A scalable and implicit meshless {RBF} method for the
  3{D} unsteady nonlinear {R}ichards equation with single and multi-zone
  domains, International journal for numerical methods in engineering 85~(2)
  (2011) 135--163.
\newblock \href {https://doi.org/https://doi.org/10.1002/nme.2960}
  {\path{doi:https://doi.org/10.1002/nme.2960}}.

\bibitem{boujoudar2022modelling}
M.~Boujoudar, A.~Beljadid, A.~Taik, Modelling of unsaturated flow through
  porous media using meshless methods, in: Canadian Society of Civil
  Engineering Annual Conference, Springer, 2022, pp. 565--576.

\bibitem{protopapas1991analytical}
A.~L. Protopapas, R.~L. Bras, Analytical solutions for unsteady
  multidimensional infiltration in heterogeneous soils, Water resources
  research 27~(6) (1991) 1029--1034.
\newblock \href {https://doi.org/https://doi.org/10.1029/91WR00331}
  {\path{doi:https://doi.org/10.1029/91WR00331}}.

\bibitem{yeh1989one}
T.-C.~J. Yeh, One-dimensional steady state infiltration in heterogeneous soils,
  Water Resources Research 25~(10) (1989) 2149--2158.
\newblock \href {https://doi.org/https://doi.org/10.1029/WR025i010p02149}
  {\path{doi:https://doi.org/10.1029/WR025i010p02149}}.

\bibitem{merrill1978laterally}
S.~Merrill, P.~Raats, C.~Dirksen, Laterally confined flow from a point source
  at the surface of an inhomogeneous soil column, Soil Science Society of
  America Journal 42~(6) (1978) 851--857.
\newblock \href
  {https://doi.org/https://doi.org/10.2136/sssaj1978.03615995004200060002x}
  {\path{doi:https://doi.org/10.2136/sssaj1978.03615995004200060002x}}.

\bibitem{tartakovsky2003unsaturated}
D.~M. Tartakovsky, Z.~Lu, A.~Guadagnini, A.~M. Tartakovsky, Unsaturated flow in
  heterogeneous soils with spatially distributed uncertain hydraulic
  parameters, Journal of Hydrology 275~(3-4) (2003) 182--193.
\newblock \href {https://doi.org/https://doi.org/10.1016/S0022-1694(03)00042-8}
  {\path{doi:https://doi.org/10.1016/S0022-1694(03)00042-8}}.

\bibitem{bakker2004two}
M.~Bakker, J.~L. Nieber, Two-dimensional steady unsaturated flow through
  embedded elliptical layers, Water Resources Research 40~(12) (2004).
\newblock \href {https://doi.org/https://doi.org/10.1029/2004WR003295}
  {\path{doi:https://doi.org/10.1029/2004WR003295}}.

\bibitem{zhang2021finite}
Z.~Zhang, W.~Wang, C.~Gong, T.-c.~J. Yeh, L.~Duan, Z.~Wang, Finite analytic
  method: Analysis of one-dimensional vertical unsaturated flow in layered
  soils, Journal of Hydrology 597 (2021) 125716.
\newblock \href {https://doi.org/https://doi.org/10.1016/j.jhydrol.2020.125716}
  {\path{doi:https://doi.org/10.1016/j.jhydrol.2020.125716}}.

\bibitem{boujoudar2021localized}
M.~Boujoudar, A.~Beljadid, A.~Taik, Localized {MQ}-{RBF} meshless techniques
  for modeling unsaturated flow, Engineering Analysis with Boundary Elements
  130 (2021) 109--123.
\newblock \href
  {https://doi.org/https://doi.org/10.1016/j.enganabound.2021.05.011}
  {\path{doi:https://doi.org/10.1016/j.enganabound.2021.05.011}}.

\bibitem{kansa1990multiquadrics}
E.~J. Kansa, Multiquadrics—a scattered data approximation scheme with
  applications to computational fluid-dynamics—ii solutions to parabolic,
  hyperbolic and elliptic partial differential equations, Computers \&
  mathematics with applications 19~(8-9) (1990) 147--161.

\bibitem{lee2003local}
C.~K. Lee, X.~Liu, S.~C. Fan, Local multiquadric approximation for solving
  boundary value problems, Computational Mechanics 30~(5-6) (2003) 396--409.

\bibitem{li2013localized}
M.~Li, W.~Chen, C.~Chen, The localized {RBF}s collocation methods for solving
  high dimensional {PDE}s, Engineering Analysis with Boundary Elements 37~(10)
  (2013) 1300--1304.
\newblock \href
  {https://doi.org/https://doi.org/10.1016/j.enganabound.2013.06.001}
  {\path{doi:https://doi.org/10.1016/j.enganabound.2013.06.001}}.

\bibitem{kansa1990multiquadrics1}
E.~J. Kansa, Multiquadrics—a scattered data approximation scheme with
  applications to computational fluid-dynamics—i surface approximations and
  partial derivative estimates, Computers \& Mathematics with applications
  19~(8-9) (1990) 127--145.

\bibitem{mirinejad2017rbf}
H.~Mirinejad, T.~Inanc, An rbf collocation method for solving optimal control
  problems, Robotics and Autonomous Systems 87 (2017) 219--225.

\bibitem{vsarler2007global}
B.~{\v{S}}arler, From global to local radial basis function collocation method
  for transport phenomena, in: Advances in meshfree techniques, Springer, 2007,
  pp. 257--282.

\bibitem{hamaidi2016space}
M.~Hamaidi, A.~Naji, A.~Charafi, Space--time localized radial basis function
  collocation method for solving parabolic and hyperbolic equations,
  Engineering Analysis with Boundary Elements 67 (2016) 152--163.
\newblock \href
  {https://doi.org/https://doi.org/10.1016/j.enganabound.2016.03.009}
  {\path{doi:https://doi.org/10.1016/j.enganabound.2016.03.009}}.

\bibitem{stevens2009order}
D.~Stevens, H.~Power, H.~Morvan, An order-n complexity meshless algorithm for
  transport-type pdes, based on local hermitian interpolation, Engineering
  analysis with boundary elements 33~(4) (2009) 425--441.

\bibitem{ben2018radial}
E.~Ben-Ahmed, M.~Sadik, M.~Wakrim, Radial basis function partition of unity
  method for modelling water flow in porous media, Computers \& Mathematics
  with Applications 75~(8) (2018) 2925--2941.

\bibitem{cueto2008nonlocal}
L.~Cueto-Felgueroso, R.~Juanes, Nonlocal interface dynamics and pattern
  formation in gravity-driven unsaturated flow through porous media, Physical
  Review Letters 101~(24) (2008) 244504.
\newblock \href
  {https://doi.org/https://doi.org/10.1103/PhysRevLett.101.244504}
  {\path{doi:https://doi.org/10.1103/PhysRevLett.101.244504}}.

\bibitem{beljadid2020continuum}
A.~Beljadid, L.~Cueto-Felgueroso, R.~Juanes, A continuum model of unstable
  infiltration in porous media endowed with an entropy function, Advances in
  Water Resources 144 (2020) 103684.
\newblock \href
  {https://doi.org/https://doi.org/10.1016/j.advwatres.2020.103684}
  {\path{doi:https://doi.org/10.1016/j.advwatres.2020.103684}}.

\bibitem{keita2021implicit}
S.~Keita, A.~Beljadid, Y.~Bourgault, Implicit and semi-implicit second-order
  time stepping methods for the {R}ichards equation, Advances in Water
  Resources 148 (2021) 103841.
\newblock \href
  {https://doi.org/https://doi.org/10.1016/j.advwatres.2020.103841}
  {\path{doi:https://doi.org/10.1016/j.advwatres.2020.103841}}.

\bibitem{leverett1941capillary}
M.~Leverett, Capillary behavior in porous solids, Transactions of the AIME
  142~(01) (1941) 152--169.
\newblock \href {https://doi.org/https://doi.org/10.2118/941152-G}
  {\path{doi:https://doi.org/10.2118/941152-G}}.

\bibitem{bentley1975multidimensional}
J.~L. Bentley, Multidimensional binary search trees used for associative
  searching, Communications of the ACM 18~(9) (1975) 509--517.

\bibitem{yao2012assessment}
G.~Yao, B.~{\v{S}}arler, et~al., Assessment of global and local meshless
  methods based on collocation with radial basis functions for parabolic
  partial differential equations in three dimensions, Engineering analysis with
  boundary elements 36~(11) (2012) 1640--1648.
\newblock \href
  {https://doi.org/https://doi.org/10.1016/j.enganabound.2012.04.012}
  {\path{doi:https://doi.org/10.1016/j.enganabound.2012.04.012}}.

\bibitem{young2016localized}
D.~Young, S.~Hu, C.~Wu, Localized radial basis function scheme for
  multidimensional transient generalized newtonian fluid dynamics and heat
  transfer, Engineering Analysis with Boundary Elements 64 (2016) 68--89.
\newblock \href
  {https://doi.org/https://doi.org/10.1016/j.enganabound.2015.11.004}
  {\path{doi:https://doi.org/10.1016/j.enganabound.2015.11.004}}.

\bibitem{li2022efficient}
M.~Li, O.~Nikan, W.~Qiu, D.~Xu, An efficient localized meshless collocation
  method for the two-dimensional burgers-type equation arising in fluid
  turbulent flows, Engineering Analysis with Boundary Elements 144 (2022)
  44--54.
\newblock \href
  {https://doi.org/https://doi.org/10.1016/j.enganabound.2022.08.007}
  {\path{doi:https://doi.org/10.1016/j.enganabound.2022.08.007}}.

\bibitem{simunek2005hydrus}
J.~Simunek, M.~T. Van~Genuchten, M.~Sejna, The {HYDRUS}-1{D} software package
  for simulating the one-dimensional movement of water, heat, and multiple
  solutes in variably-saturated media, University of California-Riverside
  Research Reports 3 (2005) 1--240.

\bibitem{baron2017adaptive}
V.~Baron, Y.~Coudi{\`e}re, P.~Sochala, Adaptive multistep time discretization
  and linearization based on a posteriori error estimates for the {R}ichards
  equation, Applied Numerical Mathematics 112 (2017) 104--125.
\newblock \href {https://doi.org/https://doi.org/10.1016/j.apnum.2016.10.005}
  {\path{doi:https://doi.org/10.1016/j.apnum.2016.10.005}}.

\bibitem{lenhard1989corey}
R.~Lenhard, J.~Parker, S.~Mishra, On the correspondence between {B}rooks
  {C}orey and van {G}enuchten models, J. Irrig. Drain Eng 115~(4) (1989)
  744--751.

\end{thebibliography}

\end{document}